\documentclass[12pt,reqno]{amsart}
\usepackage{amssymb,delarray}
\usepackage{amsfonts}
\usepackage{epsfig}
\usepackage[all]{xy}
\usepackage{amscd}

\makeindex{}

\newtheorem{thm}{Theorem}[section]
\newtheorem{lem}[thm]{Lemma}
\newtheorem{prop}[thm]{Proposition}

\newtheorem{cor}[thm]{Corollary}

\newtheorem{ex}[thm]{Example}
\newtheorem{rem}[thm]{Remark}

{\catcode`\@=11
\gdef\n@te#1#2{\leavevmode\vadjust{%
 {\setbox\z@\hbox to\z@{\strut#1}%
  \setbox\z@\hbox{\raise\dp\strutbox\box\z@}\ht\z@=\z@\dp\z@=\z@%
  #2\box\z@}}}
\gdef\leftnote#1{\n@te{\hss#1\quad}{}}
\gdef\rightnote#1{\n@te{\quad\kern-\leftskip#1\hss}{\moveright\hsize}}
\gdef\?{\FN@\qumark}
\gdef\qumark{\ifx\next"\DN@"##1"{\leftnote{\rm##1}}\else
 \DN@{\leftnote{\rm??}}\fi{\rm??}\next@}}

\begin{document}

\baselineskip=14.pt plus 2pt 

\title[Alexander modules]{Alexander modules of irreducible $C$-groups 
}
\author[Vik.S.~Kulikov]{Vik. S.~Kulikov}
\address{Steklov Mathematical Institute\\
Gubkina str., 8\\
119991 Moscow \\
Russia} \email{kulikov@mi.ras.ru}

\dedicatory{} \subjclass{}
\thanks{The work  was partially supported
by the RFBR  (05-01-00455),   NWO-RFBR
047.011.2004.026 
(RFBR 05-02-89000-$NWO_a$), and by RUM1-2692-MO-05. \\ The
research was done during the visit of the author at the Max Planck
Institute of Mathematics at Bonn. The author thanks MPIM for their
hospitality and support.}
\begin{abstract}
A complete description of the Alexander modules of knotted
$n$-manifolds in the sphere $S^{n+2}$, $n\geq 2$, and irreducible
Hurwitz curves is given. This description is applied to
investigate properties of the first homology groups of cyclic
coverings of the sphere $S^{n+2}$ and the projective complex plane
$\mathbb C\mathbb P^2$ branched respectively alone knotted
$n$-manifolds and along irreducible Hurwitz (in particular,
algebraic) curves.
\end{abstract}

\maketitle
\setcounter{tocdepth}{2}


\def\st{{\sf st}}

\section*{Introduction}

A class $\mathcal C$ of $C$-groups and its subclass $\mathcal H$
of Hurwitz $C$-groups (see definitions below) play very important
role in geometry of codimension two submanifolds. For example, it
is well known that the knot  and link groups (given by Wirtinger
presentations) are $C$-groups and any $C$-group $G$ can be
realized as the group of a linked $n$-manifold if $n\geq 2$, that
is, as the fundamental group $\pi_1(S^{n+2}\setminus V)$ of the
complement of a closed oriented manifold $V$ without boundary,
$\dim_{\mathbb R}V=n$, in the $(n+2)$-dimensional sphere $S^{n+2}$
(see \cite{Ku1}) and viceversa. Note also that a $C$-group $G$ is
isomorphic to $\pi_1(S^{n+2}\setminus S^n)$, $n\geq 3$, for some
linked $n$-dimensional spheres $S^n$ if and only if $H_2G=0$
(\cite{Ke}). Some other results related to description of groups
$\pi_1(S^{n+2}\setminus S^n)$ can be found in \cite{Le1} and
\cite{Ke-W}.

If $H\subset \mathbb C\mathbb P^2$ is an algebraic or, more
generally, Hurwitz\footnote{ The definition of Hurwitz curves can
be found in \cite{Kh-Ku} or in \cite{Gr-Ku}.} (resp.,
pseudo-holomorphic) curve of degree $m$, then the Zariski -- van
Kampen presentation of $\pi_1=\pi_1(\mathbb C\mathbb P^2\setminus
(H\cup L))$ defines on $\pi_1$ a structure of a Hurwitz $C$-group
of degree $m$, where $L$ is a line at "infinity" (that is, $L$ is
a fiber of linear projection $\text{pr}:\mathbb C\mathbb P^2\to
\mathbb C\mathbb P^1$ and it is in general position with respect
to $H$; if $H$ is a pseudo-holomorphic curve, then $\text{pr}$ is
given by a pencil of pseudo-holomorphic lines). In \cite{Ku}, it
was proved that any Hurwitz $C$-group $G$ of degree $m$ can be
realized as the fundamental group $\pi_1(\mathbb C\mathbb
P^2\setminus (H\cup L))$ for some Hurwitz (resp.
pseudo-holomorphic) curve $H$, $\deg H=2^nm$, with singularities
of the form $w^m-z^m=0$, where $n$ depends on the Hurwitz
$C$-presentation of $G$. So the class $\mathcal H$ coincides with
the class $\{ \, \pi_1(\mathbb C\mathbb P^2\setminus (H\cup L))\,
\}$ of the fundamental groups of the complements of\, "affine"
Hurwitz (resp., of\, "affine" pseudo-holomorphic) curves and it
contains the subclass of the fundamental groups of complements of
affine plane algebraic curves.

By definition, a {\it $C$-group} is a group together with  a
finite presentation
\begin{equation} \label{zeroC}
G_W=\langle x_1,\dots ,x_m \, \mid \, x_i= w_{i,j,k}^{-1}
x_jw_{i,j,k} , \, \, w_{i,j,k}\in W\, \rangle ,
\end{equation}
where $W=\{ w_{i,j,k}\in \mathbb F_m\, \mid \, 1\leq i,j\leq m,\,
\, 1\leq k\leq h(i,j)\}$ is a collection consisting  of elements
of the free group $\mathbb F_m$ generated by free generators
$x_1,\dots,x_m$  (it is possible that
$w_{i_1,j_1,k_1}=w_{i_2,j_2,k_2}$ for $(i_1,j_1,k_1)\neq
(i_2,j_2,k_2)$), and $h:\{1,\dots , m\}^2 \to \mathbb Z$ is some
function. Such a presentation is called a $C$-{\it presentation}
($C$, since all relations are conjugations). Let $\varphi_W:
\mathbb F_m\to G_W$ be the canonical epimorphism. The elements
$\varphi_{W}(x_i)\in G$, $1\leq i\leq m$, and the elements
conjugated to them are called the {\it $C$-generators} of $G$. Let
$f:G_1\to G_2$ be a homomorphism of $C$-groups. It is called a
{\it $C$-homomorphism} if the images of the $C$-generators of
$G_1$ under $f$ are $C$-generators of the $C$-group $G_2$.
$C$-groups are considered up to $C$-isomorphisms. Properties of
$C$-groups were investigated in \cite{Ku3}, \cite{Ku},
\cite{Kuz},\cite{Ku.O}.

A $C$-presentation (\ref{zeroC}) is called a {\it Hurwitz
$C$-presentation of degree} $m$ if for each $i=1,\dots,m$ the word
$w_{i,i,1}$ coincides with the product $x_1\dots x_m$, and a
$C$-group $G$ is called a {\it Hurwitz $C$-group} ({\it of degree}
$m$) if for some $m\in \mathbb N$ it possesses a Hurwitz
$C$-presentation of degree $m$. In other words, a $C$-group $G$ is
a Hurwitz $C$-group of degree $m$ if there are $C$-generators
$x_1,\dots, x_m$ generating $G$ such that the product $x_1\dots
x_m$ belongs to the center of $G$. Note that the degree of a
Hurwitz $C$-group $G$ is not defined canonically and depends on
the Hurwitz $C$-presentation of $G$. Denote by $\mathcal H$ the
class of all Hurwitz $C$-groups.

It is easy to show that $G/G'$ is a finitely generated free
abelian group for any $C$-group $G$, where $G'=[G,G]$ is the
commutator subgroup of $G$. A $C$-group $G$ is called {\it
irreducible} if $G/G'\simeq \mathbb Z$ and we say that $G$
consists of $k$ {\it irreducible components} if $G/G'\simeq
\mathbb Z^k$. If a Hurwitz $C$-group $G$ is realized as the
fundamental group $\pi_1(\mathbb C\mathbb P^2\setminus (H\cup L))$
of the complement of some Hurwitz curve $H$, then the number of
irreducible components of $G$ is equal to the number of
irreducible components of $H$. Similarly, if a $C$-group $G$
consisting of $k$ irreducible components is realized as the group
of a linked $n$-manifold $V$, $G=\pi_1(S^{n+2}\setminus V)$, then
the number of connected components of $V$ is equal to $k$.

A free group $\mathbb F_n$ with fixed free generators is a
$C$-group and for any $C$-group $G$ the canonical $C$-epimorphism
$\nu : G\to \mathbb F_1$, sending the $C$-generators of $G$ to the
$C$-generator of $\mathbb F_1$, is well defined. Denote by $N$ its
kernel. Note that if $G$ is an irreducible $C$-group, then $N$
coincides with $G'$. In what follows we consider only the
irreducible case.

Let $G$ be an irreducible $C$-group. The $C$-epimorphism $\nu$
induces the following exact sequence of groups
$$1\to G^{\prime}/G^{\prime\prime}\to G/G^{\prime\prime}
\buildrel{\nu_*}\over\longrightarrow  \mathbb F_1\to 1,$$ where
$G^{\prime\prime}=[G',G']$. The $C$-generator of $\mathbb F_1$
acts on $G^{\prime}/G^{\prime\prime}$ by conjugation $\widetilde
x^{-1} g\widetilde{x}$, where $g\in G^{\prime}$ and
$\widetilde{x}$ is one of the $C$-generators of $G$. Denote by $t$
this action. The group $A_0(G)=G^{\prime}/G^{\prime\prime}$ is an
abelian group and the action $t$ defines on $A_0(G)$ a structure
of $\Lambda$-module, where $\Lambda=\mathbb Z[t,t^{-1}]$ is the
ring of Laurent polynomials with integer coefficients. The
$\Lambda$-module $A_0(G)$ is called the {\it Alexander module} of
the $C$-group $G$. The action $t$ induces an action $h_{\mathbb
C}$ on $A_{\mathbb C}=A_0(G)\otimes \mathbb C$ and it is easy to
see that its characteristic polynomial $h_{\mathbb C}\in \mathbb
Q[t]$. The polynomial $\Delta (t)=a\det (h_{\mathbb
C}-t\text{Id})$, where $a\in \mathbb N$ is the smallest number
such that $a\det (h_{\mathbb C}-t\text{Id})\in \mathbb Z[t]$, is
called the {\it Alexander polynomial} of the $C$-group $G$. If $H$
is either an algebraic, or Hurwitz, or pseudo-holomorphic
irreducible curve in $\mathbb C\mathbb P^2$ (resp., $V\subset
S^{n+2}$ is a knotted (that is, connected smooth oriented without
boundary) $n$-manifold, $n\geq 1$) and $G=\pi_1(\mathbb C\mathbb
P^2\setminus (H\cup L))$ (resp., $G=\pi_1(S^{n+2}\setminus V)$),
then the Alexander module $A_0(G)$ of the group $G$ and its
Alexander polynomial $\Delta (t)$ are called the {\it Alexander
module and Alexander polynomial} of the curve $H$ (resp., of the
knotted manifold $V$). Note that the Alexander module $A_0(H)$ and
the Alexander polynomial $\Delta (t)$ of a curve $H$ do not depend
on the choice of the generic (pseudo)-line $L$. Results related to
the Alexander modules of knotted spheres are stated in \cite{Le2},
\cite{Le3}.

In \cite{Gr-Ku} and \cite{Ku?}, properties of the Alexander
polynomials of Hurwitz curves were investigated. In particular, it
was proved that  if $H$ is an irreducible Hurwitz curve of degree
$d$, then its Alexander polynomial $\Delta (t)$ has the following
properties
\begin{itemize}
\item[($i$)] $\Delta (t)\in \mathbb Z[t]$,  $\deg \Delta (t)$ is
an even number; \item[($ii$)] $\Delta(0)=\Delta (1)=1$;
\item[($iii$)] $\Delta (t)$ is a divisor of the polynomial
$(t^d-1)^{d-2}$,
\end{itemize}
and, moreover, a polynomial $P(t)\in \mathbb Z[t]$ is the
Alexander polynomial of an irreducible Hurwitz curve if and only
if the roots of $P(t)$ are roots of unity and $P(1)=1$. \\

Let $G=\pi_1(\mathbb C\mathbb P^2\setminus(H\cup L))$ be the
fundamental group of the complement of an irreducible affine
Hurwitz curve (resp., $G=\pi_1(S^{n+2}\setminus V)$ is the group
of a knotted $n$-manifold, $n\geq 1$). The homomorphism $\nu :G
\to \mathbb F_1$ defines an infinite unramified cyclic covering
$f_{\infty}:X_{\infty}\to \mathbb C\mathbb P^2\setminus (H\cup L)$
(resp., $f_{\infty}:X_{\infty}\to S^{n+2}\setminus V$). We have
$H_1(X_{\infty},\mathbb Z)= G^{\prime}/G^{\prime\prime}$ and the
action of $t$ on $H_1(X_{\infty},\mathbb Z)$ coincides with the
action of a generator $h$ of the covering transformation group of
the covering $f_{\infty}$.

For any $k\in \mathbb N$ denote\ by\! $\text{mod}_{k}:\mathbb
F_1\to \mu_k=\mathbb F_1/\{t^k\}$ the\ na\-tu\-ral epimorphism to
the cyclic group $\mu_k$ of degree $k$. The covering $f_{\infty}$
can be factorized through the cyclic covering $f'_k:X'_k\to
\mathbb C\mathbb P^2\setminus (H\cup L)$ (resp., $f'_k:X'_k\to
S^{n+2}\setminus V$) associated with the epimorphism\
$\text{mod}_k\circ \nu$, $f_{\infty}=f'_k\circ g_k$. Since a
Hur\-witz curve $H$ has only analytic singularities, the covering
$f'_k$ can be extended (see \cite{Gr-Ku}) to a map ${\widetilde
f}_k:{\widetilde X}_k\to X$ branched along $H$ and, maybe, along
$L$. Here $\widetilde X_k$ is a closed four dimensional variety
locally isomorphic over a singular point of $H$ to a complex
analytic singularity given by an equation $w^k=F(u,v)$, where
$F(u,v)$ is a local equation of $H$ at its singular point. In
addition, $\widetilde X_k$ is locally isomorphic over a
neighbourhood of an intersection point of $H$ and $L$ to the
singularity locally given by $w^k=vu^{d}$, where $d$ is the
smallest non-negative integer for which $m+d$ is divisible by $k$.
The variety $\widetilde X_k$, if $\widetilde f_k^{-1}(L)\subset
\text{Sing}\, \widetilde X_k$, can be normalized (as in the
algebraic case) and we obtain a covering $\widetilde f_{k,\,
\text{norm}}:\widetilde X_{k,\, \text{norm}}\to\mathbb C\mathbb
P^2$ in which $\widetilde X_{k,\, \text{norm}}$ is a singular
analytic variety at its finitely many singular points. The map
${\widetilde f}_{k,\, \text{norm}}$ is branched along $H$ and,
maybe, along the line "at infinity" $L$ (if $k$ is not a divisor
of $\deg H$, then ${\widetilde f}_{k,\, \text{norm}}$ is branched
along $L$). One can resolve the singularities of $\widetilde
X_{k,\, \text{norm}}$ and obtain a smooth manifold $\overline
X_k$, $\dim _{\mathbb R}\overline X_k=4$. Let $\sigma :\overline
X_k\to \widetilde X_{k,\, \text{norm}}$ be a resolution of the
singularities, $E=\sigma^{-1}(\text{Sing}\, \widetilde X_{k,\,
\text{norm}})$ the preimage of the set of singular points of
$\widetilde X_{k,\, \text{norm}}$, and $\overline f_k=\widetilde
f_{k,\, \text{norm}}\circ \sigma$.  The action $h$ induces an
action $\overline h_k$ on $\overline X_{k}$ and an action $t$ on
$H_1(\overline X_{k},\mathbb Z)$.

Similarly, the covering $f'_k:X'_k\to S^{n+2}\setminus V$ can be
extended to a smooth map $f_k:X_k\to S^{n+2}$ branched along $V$,
where $X_k$ is a smooth compact $(n+2)$-manifold, and the action
$t$ induces actions $h_k$ on $X_{k}$ and $h_{k*}$ on
$H_1(X_{k},\mathbb Z)$. The action $h_{k*}$ defines on
$H_1(X_{k},\mathbb Z)$ a structure of $\Lambda$-module.

In \cite{Gr-Ku}, it was shown that for any Hurwitz curve $H$, a
covering space $\overline X_k$ can be embedded as a symplectic
submanifold to a complex projective rational 3-fold on which the
symplectic structure is given by an integer K\"{a}hler form, and
it was proved that the first Betti number $b_1(\overline
X_k)=\dim_{\mathbb C} H_1(\overline X_k,\mathbb C)$ of $\overline
X_k$ is equal to $r_{k,\neq 1}$, where $r_{k,\neq 1}$ is the
number of roots of the Alexander polynomial $\Delta (t)$ of the
curve $\bar H$ which are $k$-th roots of unity not
equal to 1. \\

Let $M$ be a Noetherian $\Lambda$-module. We say that $M$ is {\it
$(t-1)$-invertible} if the multiplication by $t-1$ is an
automorphism of $M$. A $\Lambda$-module $M$ is called {\it
$t$-unipotent} if for some $n\in \mathbb N$ the multiplication by
$t^n$ is the identity automorphism of $M$.  The smallest $k\in
\mathbb N$ such that $$t^k-1\in \text{Ann}(M)=\{ f(t)\in \Lambda
\mid f(t)v=0 \, \, \text{for}\, \, \forall v\in M\}$$ is called
the {\it unipotence index} of $t$-unipotent module $M$.

Let $M$ be a Noetherian $(t-1)$-invertible $\Lambda$-module. A
$t$-invertible $\Lambda$-modules $A_n(M)=M/(t^k-1)M$ is called the
{\it $k$-th derived Alexander module} of $M$ and if $M$ is the
Alexander module of a $C$-group $G$ (resp., of a knotted
$n$-manifold $V$, resp., of a Hurwitz curve $H$), then $A_k(M)$ is
called the {\it $k$-th derived Alexander module }of $G$ (resp., of
$V$, resp., of $H$) and it will be denoted by $A_k(G)$ (resp.,
$A_k(V)$, resp., $A_k(H)$)

The main results of the article are the following statements.

\begin{thm} \label{main1} A  $\Lambda$-module $M$ is the Alexander
module of a knotted $n$-manifold, $n\geq 2$, if and only if it is
a Noetherian $(t-1)$-invertible $\Lambda$-module.
\end{thm}

\begin{thm} \label{kmain1} Let $V$ be a knotted $n$-mani\-fold,
$n\geq 1$, and $f_k:X_k\to S^{n+2}$ the cyclic covering branched
along $V$. Then $H_1(X_k,\mathbb Z)$ is isomorphic to the $k$-th
Alexander module $A_k(V)$ of $V$ as a $\Lambda$-module.
\end{thm}

Similar statements hold in the case of algebraic and, more
generally, of Hurwitz (resp., pseudo-holomorphic) curves.

\begin{thm} \label{main2} A  $\Lambda$-module $M$ is the Alexander
module of an irreducible Hurwitz {\rm (}resp.,
pseudo-holomorphic{\rm )} curve if and only if it is a Noetherian
$(t-1)$-invertible $t$-unipotent $\Lambda$-module. In particular,
the Alexander module of an irreducible algebraic plane curve is a
Noetherian $(t-1)$-invertible $t$-unipotent $\Lambda$-module.

The unipotence index of the Alexander module $A_0(H)$ of an
irreducible plane algebraic {\rm (}resp., Hurwitz or
pseudo-holomorphic{\rm )} curve $H$ is a divisor of $\deg H$.
\end{thm}

\begin{cor} \label{finger} The Alexander module $A_0(H)$ of an
irreducible plane algebraic {\rm (}resp., Hurwitz or
pseudo-holomorphic{\rm )} curve $H$ is finitely generated over
$\mathbb Z$, that is, $A_0(H)$ is a finitely generated abelian
group.

A finitely generated abelian group $G$ is the Alexander module
$A_0(H)$ of some irreducible Hurwitz or pseudo-holomorphic curve
$H$ if and only if there are an integer $m$ and an automorphism
$h\in \text{Aut}(G)$ such that $h^m=\text{Id}$ and $h-\text{Id}$
is also an automorphism of $G$.
\end{cor}

\begin{thm} \label{hmain2}
Let $H$ be an algebraic {\rm (}resp., Hurwitz or
pseudo-holomor\-phic{\rm )} irreducible curve in $\mathbb C\mathbb
P^2$, $\deg H=m$, and ${\overline f}_k:{\overline X}_k\to \mathbb
C\mathbb P^2$ be a resolution of singularities of the cyclic
covering of degree $\deg \overline f_k=k$ branched along $H$ and,
maybe, alone the line "at infinity" $L$. Then
$$\begin{array}{l} H_1(\overline X_k\setminus E,\mathbb Z)\simeq A_k(H), \\
H_1(\overline X_k,\mathbb Q)\simeq A_k(H)\otimes \mathbb Q,
\end{array}$$ where $A_k(H)$ is the $k$-th Alexander module of $H$ and
$E=\sigma^{-1}({\rm Sing}\, \widetilde X_{k,\, {\rm norm }})$.
\end{thm}

It should be noticed that in general case  the homomorphism
$H_1(\overline X_k\setminus E,\mathbb Z)\simeq A_k(H) \to
H_1(\overline X_k,\mathbb Z)$, induced by the embedding $\overline
X_k\setminus E\hookrightarrow \overline X_k$, is an epimorphism
and it is not necessary to be an isomorphism (see Example
\ref{exK3}).

\begin{cor} \label{ckmain2} Let $H$ be an algebraic
(resp, Hurwitz or pseudo-holomorphic) irreducible curve in
$\mathbb C\mathbb P^2$, $\deg H=m$, and ${\overline
f}_k:{\overline X}_k\to \mathbb C\mathbb P^2$ be a resolution of
singularities of the cyclic covering of degree $\deg f_k=k$
branched along $H$ and, maybe, alone the line "at infinity". Then
\begin{itemize}
\item[$(i)$] the first Betti number $b_1(\overline X_k)$ of
$\overline X_k$ is an even number; \item[$(ii)$] if $k=p^r$, where
$p$ is prime, then $H_1(\overline X_k,\mathbb Q)=0$;
\item[$(iii)$] if $k$ and $m$ are coprime, then $H_1(\overline
X_k,\mathbb Z)=0$; \item[$(iv)$] $H_1(\overline X_2,\mathbb Z)$ is
a finite abelian group of odd order.
\end{itemize}
\end{cor}

Note also that any $C$-group $G$ can be realized (see \cite{Ku})
as $\pi_1(\Delta^2\setminus (C\cap \Delta^2))$, where $\Delta^2=\{
|z|<1\}\times \{ |w|<1\}\subset \mathbb C^2$ is a bi-disc and
$C\subset \mathbb C^2$ is a non-singular algebraic curve such that
the restriction of $\text{pr}_1:\Delta^2\to \{ |z|<1\}$ to $C\cap
\Delta^2$ is a proper map. Therefore the analogue of Theorems
\ref{main1} and \ref{kmain1} and corollaries of them hold also in
this case.

The proof of Theorems \ref{main1} and  \ref{main2} is given in
section 3. In section 1, properties of Noetherian
$(t-1)$-invertible $\Lambda$-modules are described and section 2
is devoted to Noetherian $t$-unipotent $\Lambda$-modules. In
section 4, Theorems \ref{kmain1} and \ref{hmain2} are proved and
some other corollaries of them are stated.  \\

\section{$(t-1)$-invertible $\Lambda$-modules }
\subsection{Criteria of $(t-1)$-invertibility}
Before to describe  $(t-1)$-invertible $\Lambda$-modules, let us
recall that the ring $\Lambda=\mathbb Z[t,t^{-1}]$ is Noetherian.
Each element $f\in \Lambda$ can be written in the form
$$f=\sum_{n_{-}\leq i\leq n_{+}}a_it^i\in \mathbb Z[t,t^{-1}],$$
where $n_-,n_+,i,a_{i}\in \mathbb Z$. If $n_{-}\geq 0$ for
$f\in\Lambda$, then $f\in\mathbb Z[t]$ and it will be called  a
{\it polynomial}.

For any $n\in \mathbb Z$, $n\neq 0$, a $\mathbb Z$-homomorphism
$$f(t)=\sum
a_it^i\mapsto f(n)=\sum
a_in^i$$ is well defined. The image $f(n)$ of $f(t)$ is called the
{\it value} of $f(t)$ at $n$. If $f(t)$ is a polynomial, then its
value $f(0)=a_0$ is also well defined.

We begin with the following lemma.

\begin{lem} \label{crit0}
A Noetherian $\Lambda$-module $M$ is $(t-1)$-invertible if and
only if the multiplication by $t-1$ is a surjective endomorphism
of $M$.
\end{lem}
\proof Lemma follows from some more general statement. Namely, any
surjective $\Lambda$-endomorphism $f:M\to M$ of a Noetherian
$\Lambda$-module $M$ is an isomorphism. Indeed, if $\ker f\neq 0$,
then the chain of submodules
$$ \ker f\subset \ker f^2\subset \dots \subset \ker f^n\subset
\dots $$ is strictly increasing, since $f$ is an epimorphism. This
contradicts the Noetherian property of the module $M$. \qed \\

Let $M$ be a Noetherian $(t-1)$-invertible $\Lambda$-module.
Consider an element $v\in M$ and denote by $M_v=<v>$ a principal
submodule of $M$ generated by $v$. Since $M$ is Noetherian, any
principle submodule of $M$ is contained in a maximal principle
submodule of $M$.
\begin{lem} \label{principle}
Any maximal principal submodule $M_v$ of $(t-1)$-invertible module
$M$ is $(t-1)$-invertible.
\end{lem}
\proof Since $M$ is $(t-1)$-invertible module, there is an element
$v_1\in M$ such that $v=(t-1)v_1$. Therefore $M_v\subset M_{v_1}$.
Since $M_v$ is a maximal principle submodule of $M$, we have
$M_v=M_{v_1}$. Therefore $v_1\in M_v$ and the multiplication by
$t-1$ defines a surjective endomorphism of $M_v$. To complete the
proof we apply Lemma \ref{crit0}. \qed \\

A principal submodule $M_v\subset M$ is isomorphic to
$\Lambda/\text{Ann}_v$, where $\text{Ann}_v=\{ f\in \Lambda |
fv=0\}$ is the {\it annihilator} of $v$. The annihilator
$\text{Ann}_v$ of an element $v\in M$ is an ideal of $\Lambda$.
Denote by
$$\text{Ann}(M)=\bigcap_{v\in M}\text{Ann}_v=\{ g(t)\in
\Lambda\mid g(t)v=0\,\, \text{for}\, \, \forall v\in M\} $$ the
{\it annihilator} of $M$.

\begin{lem} \label{kriteri}
A principal $\Lambda$-module $M=\Lambda /I$ is a
$(t-1)$-invertible if and only if the ideal $I$ contains a
polynomial $f(t)$ such that $f(1)=1$.
\end{lem}
\proof Let $M$ is generated by an element $v\in M$.

If a polynomial $f(t)$ such that $f(1)=1$ is contained in
$I=\text{Ann}_v$, then $f(t)$ can be expressed in the form
\begin{equation} \label{eq1} f(t)=(t-1)g(t)+1
\end{equation}
for some polynomial $g(t)$. Therefore $v=(t-1)v_1$, where
$v_1=-g(t)v$. Thus, the multiplication by $t-1$ is a surjective
automorphism of $M$ and hence, by Lemma \ref{crit0}, the
multiplication by $t-1$ is an isomorphism of $M$.

Conversely, if $M$ is $(t-1)$-invertible, then there is an element
$v_1\in M$ such that $v=(t-1)v_1$. Let $v_1=h(t)v$ for some
$h(t)\in \Lambda$. We have $(1-(t-1)h(t))v=0$. Therefore
$1-(t-1)h(t) \in \text{Ann}_v=I$. There is an integer $k$ such
that $f(t)=t^k(1-(t-1)h(t))\in I\cap \mathbb Z[t]$. It is easy to
see that $f(1)=1$. \qed \\

As a consequence of Lemma \ref{kriteri} we obtain the following
Lemma.
\begin{lem} \label{cor-krit}
Any principal submodule of a principal $(t-1)$-invertible module
$M$ is $(t-1)$-invertible.
\end{lem}
\proof Indeed, let $M$ be generated by an element $v\in M$ and its
submodule $M_1$ be generated by $v_1=h(t)v$. Then
$\text{Ann}_v\subset \text{Ann}_{v_1}$.

Since $M$ is $(t-1)$-invertible, by Lemma \ref{kriteri}, there is
a polynomial $f(t)\in \text{Ann}_v$ such that $f(1)=1$. Applying
again Lemma \ref{kriteri}, we have that $M_1$ is
$(t-1)$-invertible, since $f(t)\in \text{Ann}_{v_1}$. \qed

\begin{prop} \label{cor-krit2}
Any submodule of a Noetherian $(t-1)$-invertible $\Lambda$-module
$M$ is $(t-1)$-invertible.
\end{prop}
\proof Let $N$ is a submodule of $M$. Since $M$ is a Noetherian
$\Lambda$-module, the submodule $N$ is generated by a finite set
of elements, say $v_1,\dots,v_n$. By Lemma \ref{cor-krit}, each
principal submodule $M_{v_i}\subset N \subset M$ is
$(t-1)$-invertible. 
Therefore the multiplication by $t-1$ is a surjective endomorphism
of $N$, since it is surjective on each $M_{v_i}\subset N$ and the
elements $v_1,\dots,v_n$ generate the module $N$. To complete the
proof, we apply Lemma \ref{crit0}. \qed

\begin{prop} \label{cor-krit3}
Any factor module of a Noetherian $(t-1)$-invertible
$\Lambda$-module $M$ is $(t-1)$-invertible.
\end{prop}
\proof It follows from Lemma \ref{crit0}. \qed \\
\begin{lem} \label{sum} Let $M_1,\dots, M_k$ be Noetherian $(t-1)$-invertible
$\Lambda$-modules. Then the direct sum $M=\bigoplus_{i=1}^kM_i$ is
a Noetherian $(t-1)$-invertible $\Lambda$-module.
\end{lem}
\proof Obvious.

\begin{cor} \label{cor-krit6} Any Noetherian $(t-1)$-invertible
$\Lambda$-module $M$ is a the factor module of a direct sum
$\bigoplus_{j=1}^n \Lambda /I_j$ of principle $(t-1)$-invertible
$\Lambda$-modules $\Lambda /I_j$.
\end{cor}
\proof Since $M$ is a Noetherian $\Lambda$-module, it is generated
by a finite set of elements, say $v_1,\dots,v_n$. By Proposition
\ref{cor-krit2}, each principal submodule $M_{v_i}\subset M$ is
$(t-1)$-invertible and, obviously, there is an epimorphism
$\bigoplus_{j=1}^n M_{v_i} \mapsto M$. \qed

\begin{rem} \label{rem} An abelian group $G$ admits
a  structure of $(t-1)$-invertible $\Lambda$-module if and only if
it has an automorphism $t$ such that $t-1$ is also an
automorphism. If $G$ is finitely generated and $t\in \text{Aut}\,
G$ is chosen, then $G$ is a Noetherian $\Lambda$-module.
\end{rem}

Note that in general case an abelian group admits many structures
of $(t-1)$-invertible $\Lambda$-modules. For example, the group
$\mathbb Z/9\mathbb Z$ admits $3$ such structures: either $tv=2v$,
or $tv=5v$, or $tv=8v$, where $v$ is a generator of $\mathbb
Z/9\mathbb Z$.

\begin{thm} \label{cor-krit4}
A Noetherian $\Lambda$-module $M$ is $(t-1)$-invertible if and
only if there is a polynomial  $f(t)\in\text{Ann}(M)$ such that
$f(1)=1$.
\end{thm}
\proof If $M$ is $(t-1)$-invertible, then, by Proposition
\ref{cor-krit2}, its each principal submodule $M_v$ is also
$(t-1)$-invertible. Therefore, by Lemma \ref{kriteri}, the
annihilator $\text{Ann}_v$ of $v\in M$ contains a polynomial
$f_v(t)$ such that $f_v(1)=1$. If $M$ is generated by $v_1,\dots,
v_n$, then the polynomial $f(t)=f_{v_1}(t)\dots f_{v_n}(t)$ is a
desired one.

Let us show that if  there is a polynomial $f(t)\in\text{Ann}(M)$
such that $f(1)=1$, then $M$ is a $(t-1)$-invertible module.
Indeed, in this case by Lemma \ref{kriteri}, each principle
submodule $M_{v}$ of $M$ is $(t-1)$-invertible. Therefore the
multiplication by $t-1$ is an isomorphism of $M$, since it is an
isomorphism of each principle submodule $M_v$ of $M$. \qed \\

As a consequence of Theorem \ref{cor-krit4} we obtain that {\it
any Noetherian $(t-1)$-invertible module $M$ is a torsion
$\Lambda$-module} and, consequently, $$\dim_{\mathbb Q} M\otimes
\mathbb Q<\infty .$$

The following proposition will be used in the proof of Theorems
\ref{main1} and \ref{main2}.
\begin{prop} \label{usef} Any Noetherian $(t-1)$-invertible
$\Lambda$-module $M$ is isomorphic to a
factor module $\Lambda^n/M_1$ of a free $\Lambda$-module
$\Lambda^n$, where the submodule $M_1$ is generated by elements
$w_1,\dots,w_n,\dots,w_{n+k}$ of $\Lambda^n$ such that
\begin{itemize}
\item[($i$)] for  $i=1,\dots, n$ the vector $w_i= (0, \dots, 0,
f_i(t),0,\dots,0)$, where a polynomial $f_i(t)$ stands on the
$i$-th place and it is such that $f_i(1)=1$, \item[($ii$)]
$w_{n+j}=(t-1)\overline w_{n+j}=((t-1)g_{j,1}(t),\dots
,(t-1)g_{j,n}(t))$ for $j=1,\dots,k$, where $g_{j,l}(t)$ are
polynomials, \item[($iii$)] if for some $m\in \mathbb N$ the
polynomial $t^m-1\in \text{Ann}(M)$, then for $i=1,\dots, n$ the
vector $w_{n+i}= (0, \dots, 0,t^m-1,0,\dots,0)$, where the
polynomial $t^m-1$ stands on the $i$-th place.
\end{itemize}
\end{prop}
\proof Let us choose generators $v_1,\dots,v_n$ of the Noetherian
$\Lambda$-module $M$. Then, by Theorem \ref{cor-krit4}, there are
polynomials $f_i(t)\in \text{Ann}_{v_i}$ such that $f_i(1)=1$.
Obviously, there is an epimorphism $$h_1: \bigoplus_{i=1}^n\Lambda
/(f_i(t))\to M$$ of $\Lambda$-modules such that $h(u_i)=v_i$ for
$u_i=(0,\dots,0,1,0\dots,0)$ where $1$ stands on the $i$-th place.
The kernel $N=\ker h$ is a Noetherian $\Lambda$-module. Let it be
generated by $$u_{n+1}=(g_{1,1}(t), \dots ,g_{1,n}(t)),\, \,
\dots, \, \, u_{n+k}=(g_{k,1}(t),\dots ,g_{k,n}(t)).$$ Without
loss of generality, we can assume that all $g_{i,j}(t)$ are
polynomials.

By Theorem \ref{cor-krit4}, the $\Lambda$-module
$\bigoplus_{i=1}^n\Lambda /(f_i(t))$ is $(t-1)$-invertible and by
Proposition \ref{cor-krit2}, $N$ is also $(t-1)$-invertible
$\Lambda$-module. Therefore the elements $(t-1)
u_{n+1},\dots,(t-1) u_{n+k}$ are also generate $N$.

If for some $m\in \mathbb N$ the polynomial $t^m-1\in
\text{Ann}(M)$, then the elements $(0, \dots,
0,t^m-1,0,\dots,0)\in N$, where the polynomial $t^m-1$ stands on
the $i$-th place. Therefore we can add the elements $(0, \dots,
0,t^m-1,0,\dots,0)$ to the set $(t-1)u_{n+1},\dots,(t-1)u_{n+k}$
and renumber the elements $\overline u_{n+1},\dots, \overline
u_{n+k}$ (here we put $k:=n+k$) of the obtained set generating $N$
so that $\overline u_{n+j}=(0, \dots, 0,t^n-1,0,\dots,0)\in N$ for
$j=1,\dots,n$, where the polynomial $t^m-1$ is stands on the
$j$-th place.

Now, to complete the proof, notice that the kernel $M_1$ of the
composite map $h\circ \nu :\Lambda^n\to M$ of $h$ and the natural
epimorphism $\nu :\Lambda^n\to \bigoplus_{i=1}^n\Lambda /(f_i(t))$
is generated by the elements $$w_i= (0, \dots, 0,
f_i(t),0,\dots,0), \qquad i=1,\dots,n,$$ where the polynomial
$f_i(t)$ stands on the $i$-th place, and the elements
$$w_{n+i}=(f_{i,1}(t),\dots,f_{i,n}(t))\in \Lambda^n, \qquad
i=1,\dots,k,$$ where the coordinates $f_{i,j}(t)$ of each
$w_{n+i}$ coincide with the coordinates $\overline g_{{i},j}(t)$
of $\overline u_{n+i}=(\overline g_{i,j}(t), \dots,\overline
g_{i,j}(t))$. \qed
\subsection{$\mathbb Z$-torsion submodules of $(t-1)$-invertible $\Lambda$-modules}
An element $v$ of a $\Lambda$-module $M$ is said to be of a {\it
finite order} if there is $m\in \mathbb Z\setminus \{ 0\}$ such
that $mv=0$. A $\Lambda$-module $M$ is called {\it $\mathbb
Z$-torsion} if all its elements are of finite order. For any
$\Lambda$-module $M$ denote by $M_{fin}$ a subset of $M$
consisting of all elements of finite order. It is easy to see that
$M_{fin}$ is a $\mathbb Z$-torsion $\Lambda$-module. If $M$ is a
Noetherian $(t-1)$-invertible $\Lambda$-module, then $M_{fin}$ is
also a Noetherian $(t-1)$-invertible $\Lambda$-module, and it
follows from Propositions \ref{cor-krit2} and \ref{cor-krit3} that
there is an exact sequence of $\Lambda$-modules
$$0 \to M_{fin}\to M\to M_1\to 0$$
in which $M_1$ is a Noetherian $(t-1)$-invertible $\Lambda$-module
free from elements of finite order.

Let $M=M_{fin}$ be a Noetherian $(t-1)$-invertible
$\Lambda$-module. Since $M$ is finitely generated over $\Lambda$,
there is an integer $d\in \mathbb N$ such that $dv=0$ for all
$v\in M$ (such $d$ will be called an {\it exponent} for $M$). Let
$d=p_1^{r_1}\dots p_n^{r_n}$ be its prime factorization. Denote by
$M(p_i)$ the subset of $M$ consisting of all elements $v\in M$
such that $p_i^rv=0$ for some $r\in \mathbb N$. It is easy to see
that $M(p)$ is a $\Lambda$-submodule of $M$ and we call it the
{\it $p$-submodule of $M$.
\begin{thm} \label{decomposition}
Let $M=M_{fin}$ be a Noetherian $(t-1)$-invertible
$\Lambda$-module and $d=p_1^{r_1}\dots p_n^{r_n}$ its exponent.
Then $M$ is the direct sum
$$ M=\bigoplus_{i=1}^nM(p_i)$$ of its $p$-submodules.
\end{thm}
\proof It coincides with the proof of similar Theorem for abelian
groups (see, for example, Theorem 8.1 in \cite{L}). \qed \\

Since the ring $\Lambda=\mathbb Z[t,t^{-1}]$ is Noetherian, any
its ideal $I$ is finitely generated. Denote by
$I_{\text{pol}}=I\cap \mathbb Z[t]$ the ideal of the ring $\mathbb
Z[t]$. It is well known that $I=\Lambda I_{\text{pol}}$, that is,
any ideal $I$ of $\Lambda$ is generated by polynomials.

Recall that $\mathbb Z[t]$ is a factorial ring. Its units are
precisely the units of $\mathbb Z$, and its prime elements are
either primes of $\mathbb Z$ or polynomials $q(t)=\sum a_it^i$
which are irreducible in $\mathbb Q[t]$ and have content $1$ (that
is, the greatest common divisors of the coefficients $a_i$ of
$q(t)$ are equal to $1$). It follows from Euclidean algorithm that
for any two polynomials $q_1(t), q_2(t)\in \mathbb Z[t]$ there are
polynomials $h_1(t), h_2(t), r(t)\in \mathbb Z[t]$ and a constant
$d\in \mathbb Z$, $d\neq 0$, such that
\begin{equation} \label{Euc} h_1(t)q_1(t)+h_2(t)q_2(t)=dr(t),
\end{equation}
where $r(t)$ is the greatest common divisor of the polynomials
$q_1(t)$ and $q_2(t)$.

\begin{lem} \label{p^r}
Let $M$ be a Noetherian $(t-1)$-invertible $\Lambda$-module and
let $t^n-1\in\text{Ann}(M)$ for some $n=p^r$, where $p$ is prime.
Then $M$ is $\mathbb Z$-torsion.
\end{lem}
\proof If
$t^n-1=(t-1)(t^{n-1}+\dots +t+1)$
 belongs to $\text{Ann}(M)$, then the polynomial
$g_{n}(t)=t^{n-1}+\dots +t+1\in\text{Ann}(M)$, since $M$ is
$(t-1)$-invertible. For $n=p^r$ in the factorization
$$g_{p^r}(t)=\prod_{i=1}^{r}\Phi_{p^i}(t)=\prod_{i=1}^{r}\sum_{j=0}^{p-1} t^{jp^{i-1}}$$
each factor is an irreducible element of $\Lambda$.

By Theorem \ref{cor-krit4}, there is a polynomial
$f(t)\in\text{Ann}(M)$ such that $f(1)=1$ and if $n=p^r$ for some
prime $p$, then $f(t)$ and $g_{p^r}(t)$ have not common
irreducible divisors. Indeed, if $g(t)$ is a divisor of $f(t)$,
then we should have $g(1)=\pm 1$, since $f(1)=1$, but
$\Phi_{p^i}(1)=p$ for each $i$. Therefore, there are polynomials
$h_1(t)$, $h_2(t)$, and a constant $d\in \mathbb N$ such that
$h_1(t)f(t)+h_2(t)g_{p^r}(t)=d$ and hence if
$g_{p^r}(t)\in\text{Ann}(M)$, then $d\in\text{Ann}(M)$, that is,
$M$
is $\mathbb Z$-torsion. \qed \\
\subsection{Principle $(t-1)$-invertible $\Lambda$-modules}
Let $I$ be a non-zero ideal of the ring $\Lambda$. Denote by $I_m$
the subset of $I_{\text{pol}}$ consisting of all polynomials
$f(t)$ having the smallest degree (let $m$ be this smallest
degree). Note that if $f(t)\in I_m\setminus \{ 0\}$, then
$f(0)\neq 0$.

Consider any two polynomials $f_1(t),f_2(t)\in I_m$ and write them
in the form $f_i(t)=d_iq_i(t)$, where $d_i\in \mathbb Z$ and the
polynomials $q_i(t)$ have content $1$. We have $q_1(t)=q_2(t)$.
Indeed, for their common greatest divisor $r(t)$ we have $\deg
r(t)\leq m$ and, moreover, $\deg r(t)=m$ if and only if
$q_1(t)=q_2(t)$. On the other hand, it follows from (\ref{Euc})
that $d_2h_1(t)f_1(t)+d_1h_2(t)f_2(t)=d_1d_2dr(t)$ for some
polynomials $h_1(t),h_2(t)$. Therefore $d_1d_2dr(t)\in
I_{\text{pol}}$ and we should have $\deg r(t)=m$.

Applying again Euclidean algorithm for integers, we obtain that if
two polynomials $f_i(t)=d_iq(t)$ belong to $I_m$, then $d_0q(t)$
belongs also to $I_m$, where $d_0$ is the greatest common divisor
of $d_1$ and $d_2$. Thus there is a polynomial $f_m(t)=d_mq(t)\in
I_m$ such that any polynomial $f(t)\in I_m$ is divided by
$f_m(t)$. The polynomial $f_m(t)$ is defined uniquely up to
multiplication by $\pm 1$ and it will be called a {\it leading
generator} of $I$.

Let $I$ be a non-zero ideal of $\Lambda$ and $f(t)=d_mq(t)$ be its
leading generator. Then any polynomial $h(t)\in I$ should be
divisible by $q(t)$. Indeed, as above it is easy to show that if
$r(t)$ is the greatest common divisor of $f(t)$ and $h(t)$, then
there is a constant $d$ such that $dr(t)\in I$ and  since $\deg
q(t)$ is minimal for polynomials belonging to $I$, we should have
the equality $r(t)=q(t)$.

The above considerations give rise to the following proposition.

\begin{prop} \label{descript} Let $M=M_v$ be a principle $(t-1)$-invertible
$\Lambda$-module generated by an element $v$. Then the annihilator
$\text{Ann}_v$ is generated by a finite set of polynomials
$f_1(t),\dots, f_k(t)$, where $f_i(t)=d_iq_i(t)$,  $d_i\in \mathbb
Z$, $d_i\neq 0$, and $q_i(t)$ have content $1$ for all $i$, such
that $f_1(t),\dots, f_k(t)$ satisfy the following properties:
\begin{itemize}
 \item[($i$)] $\deg f_1< \deg f_2\leq
\dots \leq \deg f_k$, \item[($ii$)] $f_i(0)\neq 0$ for all $i$,
\item[($iii$)] $q_1(1)=1$, \item[($iv$)]  $q_1(t) \mid q_{i}(t)$
for $i=2,\dots, k$, \item[($v$)]  $|d_i|>1$ for $i=1,\dots,k-1$,
$d_k=1$, and $q_k(1)=1$.
\end{itemize}
 \end{prop}

A set of generators of $\text{Ann}_v$ is said to be {\it good} if
it satisfies properties ($i$) -- ($v$) from Proposition
\ref{descript}. We will distinguish the principal
$(t-1)$-invertible $\Lambda$-modules $M=M_v$ as follows. We say
that $M_v$ is of {\it finite type} if in a good system
$f_1(t),\dots, f_k(t)$ of generators of $\text{Ann}_v$ the leading
generator $f_1(t)\equiv d_1$ is a constant (that is, $q_1(t)\equiv
1$). A principle $\Lambda$-module $M_v$ is said to be of {\it
mixed type} if in a good system $f_1(t),\dots, f_k(t)$ of
generators of $\text{Ann}_v$ the degree of the leading generator
$f_1=d_1q_1(t)$ is greater than one and $\mid d_1\mid\geq 2$. It
follows from the above considerations that if a principle
$(t-1)$-invertible $\Lambda$-module $M=M_v$ is not of finite or
mixed types, then for the leading generator $f_1(t)=q_1(t)$ of a
good system of generators of $\text{Ann}_v$ we should have
$q_1(1)=1$ and therefore $\text{Ann}_v$ is a principle ideal
generated by $q_1(t)$, since any polynomial $h(t)\in \text{Ann}_v$
is divisible by $q_1(t)$. Such principle $(t-1)$-invertible
$\Lambda$-modules will be called {\it bi-principle}.

It is easy to see that if $M=M_v$ is a principle $\Lambda$-module
of finite type and $d_1\in \mathbb Z$ is the leading generator of
$\text{Ann}_v$, then all elements of $M$ have order $d_1$, that
is, a principle $\Lambda$-module $M_v$ is of finite type if and
only if it is $\mathbb Z$-torsion.

If $M=M_v$ is a bi-principle $\Lambda$-module, then $M$ has not
non-zero elements of finite order. Indeed, let $q(t)$ be a
generator of $\text{Ann}_v$. If an element $v_1=h(t)v$ has order
$m$, then $mh(t)\in \text{Ann}_v$, that is, $mh(t)$ is divisible
by $q(t)$. Since $t$ is a unite of $\Lambda$, we can assume that
$h(t)$ is a polynomial, and since $q(1)=1$, the polynomial $h(t)$
should be divisible by $q(t)$, that is, $v_1=0$.

If $M=M_v$ is a $\Lambda$-module of mixed type, then there is an
exact sequence of $\Lambda$-modules
$$0 \to M_1\to M\to M_2\to 0$$
in which $M_1$ is a principle $\Lambda$-module of finite type and
$M_2$ is a bi-principle $\Lambda$-module. Indeed, let $d_1q_1(t)$
be the leading generator of $\text{Ann}_v$. Put $v_1=q_1(t)v$.
Then it is easy to see that the $\Lambda$-module
$M_1=M_{v_1}\subset M$, generated by $v_1$, is of finite type and
the $\Lambda$-module $M_2=M/M_1\simeq \Lambda/(q_1)$ is
bi-principle.
\subsection{Finitely $\mathbb Z$-generated $(t-1)$-invertible
$\Lambda$-modules} Each $\Lambda$-module $M$ can be considered as
a $\mathbb Z$-module, that is as an abelian group.

\begin{prop} \label{finitenes}
A Noetherian $(t-1)$-invertible $\Lambda$-module $M$ is finitely
generated over $\mathbb Z$ if and only if there is a polynomial
$$q(t)=\sum_{i=0}^n a_it^i\in \text{Ann}(M)$$ such that $a_n=a_0=1$.
\end{prop}
\proof In the beginning, we prove Proposition \ref{finitenes} in
the case when $M=M_v$ is a principal $\Lambda$-module.

It is easy to see that if there is a polynomial $q(t)=\sum_{i=0}^n
a_it^i\in \text{Ann}_v$ such that $a_n=a_0=1$, then $M$ is
generated over $\mathbb Z$ by the elements $v,tv,\dots, t^{n-1}v$.

Let a $\Lambda$-module $M=M_v$ be finitely generated over $\mathbb
Z$ and $h_1(t)v$,$\dots$, $h_m(t)v$ its generators. Since the
multiplication by $t$ is an isomorphism of $M$, we can assume that
$h_i(t)$, $i=1,\dots, m$, are polynomials such that $h_i(0)=0$.
Put $n-1=\max (\deg h_1(t), \dots, \deg h_m(t))$. Since $h_1(t)v$,
 $\dots$, $h_m(t)v$ generate $M$ over $\mathbb Z$, there are
integers $b_1,\dots,b_m$ and $c_1,\dots,c_m$ such that
$$ v=\sum b_ih_i(t)v\quad \text{and}\quad t^nv=\sum c_ih_i(t)v.$$
Therefore the polynomials $1-\sum b_ih_i(t)$ and $t^n-\sum
c_ih_i(t)$ belong to $\text{Ann}_v$. Then the polynomial
$t^n+1-\sum (b_i+c_i)h_i(t)$ is a desired one.

In general case,  a Noetherian $(t-1)$-invertible $\Lambda$-module
$M$ is generated by a finite set of elements $v_1,\dots,v_m$, and
$M$ is finitely generated over $\mathbb Z$ if and only if for all
$v_i$ the principal submodules $M_{v_i}\subset M$ are finitely
generated over $\mathbb Z$.

If $g(t)\in \text{Ann}(M)$, then $g(t)\in \text{Ann}_{v_i}$ for
$i=1,\dots,m$. In particular, if there is $q(t)=\sum_{i=0}^n
a_it^i\in \text{Ann}(M)$ such that $a_n=a_0=1$, then all $M_{v_i}$
(and consequently, $M$) are finitely generated over $\mathbb Z$.

If for all $i$ the principal submodules $M_{v_i}\subset M$ are
finitely generated over $\mathbb Z$, then there are polynomials
$q_i(t)=\sum_{j=0}^{n_i} a_{i,j}t^j\in \text{Ann}_{v_i}$ such that
$a_{i,n_i}=a_{i,0}=1$. Put $n=\sum n_i $. Then the polynomial
$$q(t)=q_1(t)\dots q_n(t)=t^{n}+1+\sum_{j=1}^{n-1}a_jt^j\in \text{Ann}(M),$$
 since $q(t)\in\text{Ann}_{v_i}$ for all $v_i$. \qed \\

It follows from Proposition \ref{finitenes} that there are a lot
of $(t-1)$-invertible bi-principle modules $M=\Lambda /I$  which
are not finitely generated over $\mathbb Z$. More precisely, it is
easy to see that {\it a bi-principle $(t-1)$-invertible module
$M=\Lambda /I$ is finitely generated over $\mathbb Z$ if and only
if the ideal $I=\langle q(t)\rangle$ is generated by a polynomial
$q(t)=\sum_{i=0}^n a_it^i$ such that $q(1)=1$ and its coefficients
$a_0$ and $a_n$ are equal to $\pm 1$}.

For example, for each $m\in \mathbb N$ a $(t-1)$-invertible
bi-principle module
$$M_m=\Lambda/\langle (m+1)t-m\rangle$$
is not finitely generated over $\mathbb Z$.

\begin{thm} \label{fin=fin}
Let $M$ be a Noetherian $\mathbb Z$-torsion $(t-1)$-invertible
module. Then $M$ is finitely generated over $\mathbb Z$.
\end{thm}
\proof By Theorem \ref{decomposition}, $M$ is isomorphic the
direct sum $\bigoplus M(p_i)$ of a finite number of its
$p$-submodules. Therefore it suffices to prove Theorem in the case
when $M$ has exponent $p^r$, where $p$ is a prime number. Next, by
Corollary \ref{cor-krit6}, $M$ is a factor module of the direct
sum $\bigoplus_{j=1}^n \Lambda /I_j$ of principle
$(t-1)$-invertible $\Lambda$-modules $\Lambda /I_j$ and in our
case we can assume without loss of generality that each ideal
$I_j$ contains $p^{r_j}$ for some $r_j$. Thus it suffices to prove
Theorem in the case when $M=M_v$ is a principle $(t-1)$-invertible
$\Lambda$-module of exponent $p^r$, that is, $I=\text{Ann}_v$
contains a number $p^r$ and a polynomial $g(t)$ such that
$g(1)=1$.

Let $r=1$ and $g(t)=\sum a_it^i$. Denote by $g_1(t)=\sum_{p\mid
a_i} a_it^i$ and put $\overline g(t)= g(t)-g_1(t)$. Then
$\overline g(t)\in \text{Ann}_v$, since  $g(t),g_1(t)\in
\text{Ann}_v$. It is easy to see that $\overline g_1(1)$ and $p$
are coprime, since $g(1)=1$ and $g_1(1)\equiv 0 \mod p$. Moreover,
by construction, each coefficient of the polynomial $\overline
g(t)$ and $p$ are coprime.  Multiplying by $t^{-k}$, we can assume
that $\overline g(0)\neq 0$. Let $\overline g(t)=\sum_{i=0}^m
\overline a_it^i$. Since $\overline a_m$ and $p$ are coprime, one
can find integers $b_1$ and $c_1$ such that $b_1\overline
a_m+c_1p=1$. Similarly, there are integers $b_2$ and $c_2$ such
that $b_2\overline a_0+c_2p=1$. Therefore the polynomial
$(b_1t+b_2)\overline g(t)+p(c_1t^{m+1}+c_2)\in I$ and it is equal
to $h(t)=t^{m+1}+1+\sum_{i=1}^m(b_1\overline a_{i-1}+b_2\overline
a_i)t^i$. Therefore, by Proposition \ref{finitenes}, $M_v$ is
finitely generated over $\mathbb Z$.

Now consider general case of a principle $(t-1)$-invertible
$\Lambda$-module of exponent $p^r$. Assume that for any principle
$(t-1)$-invertible $\Lambda$-module $M'$ of exponent $p^{r_1}$,
where $r_1<r$, $M'$ is finitely generated over $\mathbb Z$. Let
$M=M_v$ is a principle $(t-1)$-invertible $\Lambda$-module $M$ of
exponent $p^r$. Then the submodule $M_{v_1}$ of $M$ generated by
$v_1=p^{r-1}v$ is of exponent $p$ and the factor module
$M_{\overline v}=M_v/M_{v_1}$ is of exponent $p^{r-1}$. Now, the
proof follows from the exact sequence
$$0\to M_{v_1}\to M\to M/M{v_1}\to 0. \qed$$
\begin{cor} \label{finit}
Any Noetherian $\mathbb Z$-torsion $(t-1)$-invertible module is
finite, that is, it is a finite abelian group.
\end{cor}

\begin{lem} \label{impossible}
A group $G=\bigoplus_{i=1}^n (\mathbb Z/2^{r_i}\mathbb Z)^{m_i}$
does not admit a structure of $(t-1)$-invertible $\Lambda$-module
if $r_i\neq r_j$ for $i\neq j$ and one of $m_i=1$.
\end{lem}
\proof Assume that $G$ has a structure of $(t-1)$-invertible
$\Lambda$-module. Then for any $r$ the subgroup $2^rG$ of $G$ is
its $\Lambda$-submodule and, by Propositions \ref{cor-krit2} and
\ref{cor-krit3}, $2^rG$ and $G/2^rG$ are $(t-1)$-invertible
$\Lambda$-modules. Therefore, without loss of generality, we can
assume that
$$G=(\mathbb Z/2\mathbb Z)\oplus(\bigoplus_{i=1}^n
(\mathbb Z/2^{r_i}\mathbb Z)^{m_i}),$$ where all $r_i\geq 2$ and
$m_i\geq 2$. Let us choose generators $v_1, \dots, v_{m+1}$ of
$G$, $m=\sum_{i=1}^nm_i$, so that
$$G\simeq (\mathbb Z/2\mathbb Z)v_1\oplus(\bigoplus_{i=2}^{m+1}
(\mathbb Z/2^{\overline r_i}\mathbb Z))v_i,$$ where all $\overline
r_i\geq 2$. Consider the $\mathbb Z$-submodule $\overline G$ of
$G$ consisting of all elements $v\in G$ of order $\leq 4$.
Obviously $\overline G$ is a $\Lambda$-submodule of $G$ and it is
generated over $\mathbb Z$ (and therefore over $\Lambda$) by
$\overline v_1=v_1$ and $\overline v_i=2^{\overline r_i-2}v_i$,
$i=2\dots, m+1$. It is easy to see that as an abelian group
$\overline G$ is isomorphic to
$$\overline G\simeq (\mathbb Z/2\mathbb Z)\overline v_1\oplus(\bigoplus_{i=2}^{m+1}
(\mathbb Z/4\mathbb Z))\overline v_i.$$

By Proposition \ref{cor-krit2}, $\overline G$ is
$(t-1)$-invertible $\Lambda$-module. The multiplication by $t$ is
an automorphism of $\overline G$. Let
\begin{equation} \label{equ}
\begin{array}{ll}
 t\overline v_1= & a_1\overline v_1+ \displaystyle 2\sum_{i=2}^{m+1} b_i\overline v_i, \\
t\overline v_j= & a_j\overline v_1+ \displaystyle \sum_{i=2}^{m+1}
c_{j,i}\overline v_i, \qquad j=2,\dots, m+1, \\
\end{array}
\end{equation}
where each $a_j=0$ or $1$.

Let us show that $a_1=1$. Indeed, assume that $a_1=0$. Since the
multiplication by $t$ is an automorphism and $\overline v_1,
\dots, \overline v_{m+1}$ generate $\overline G$, we should have
an equality $\overline v_1=\sum d_it\overline v_i$, where one of
$d_i$ is odd for some $i\geq 2$ if $a_1=0$. Next, the element
$\overline v_1$ is of second order, therefore $2\sum_{i=2}^{m+1}
d_it\overline v_i=0$. On the other hand, $t\overline v_2, \dots,
t\overline v_{m+1}$ are linear independent over $\mathbb
Z/4\mathbb Z$, since $\overline v_2, \dots, \overline v_{m+1}$ are
linear independent over $\mathbb Z/4\mathbb Z$ and the
multiplication by $t$ is an isomorphism. Therefore the equality
$2\sum_{i=2}^{m+1} d_it\overline v_i=0$ is impossible if some of
$d_i$ is odd, and hence $a_1$ in (\ref{equ}) should be equal to 1.

Let us show that $\overline G$ can not be $(t-1)$-invertible.
Indeed, we have $$ t\overline v_1=\overline v_1+ \displaystyle
2\sum_{i=2}^{m+1} b_i\overline v_i.$$ Therefore
$$ (t-1)\overline v_1=\displaystyle
2\sum_{i=2}^{m+1} b_i\overline v_i $$ and the above arguments show
that the multiplication by $t-1$ is not an automorphism of
$\overline G$, since $(t-1)\overline v_1$ is a linear combination
of the elements $\overline v_2,\dots,\overline v_{m+1}$. \qed

\begin{thm} \label{descript1}
An abelian group $$G=G_1\oplus(\bigoplus_{i=1}^n (\mathbb
Z/2^{r_i}\mathbb Z)^{m_i}),$$ where $r_i\neq r_j$ for $i\neq j$
and $G_1$ is a group of odd order, admits a structure of
$(t-1)$-invertible $\Lambda$-module if and only if all $m_i\geq
2$.
\end{thm}
\proof By Theorem \ref{decomposition}, if $M=M_{fin}$ is a
Noetherian $(t-1)$-invertible $\Lambda$-module and
$d=p_1^{r_1}\dots p_n^{r_n}$ its exponent, then $M$ is the direct
sum $$ M=\bigoplus_{i=1}^nM(p_i)$$ of its $p$-submodules which are
$(t-1)$-invertible by Proposition \ref{cor-krit2}. Now,  each its
submodule $M(p_i)$ with odd $p_i$ is of odd order and, by Lemma
\ref{impossible}, its $2$-submodule $M(2)$ is isomorphic (as an
abelian group) to $\bigoplus_{i=1}^k (\mathbb Z/2^{r_i}\mathbb
Z)^{m_i}$, where all $m_i\geq 2$.

To prove the inverse statement, note, first,  that the finite
direct sum of $(t-1)$-invertible $\Lambda$-modules is also a
$(t-1)$-invertible $\Lambda$-module. Next, for any prime $p>2$, a
$(t-1)$-invertible $\Lambda$-module $M=\Lambda /I$, where $I$ is
generated by the number $p^r$ and polynomial $2t-1$, is isomorphic
to $\mathbb Z/p^r\mathbb Z$ as an abelian group. Finally, for
$n\geq 2$ the $(t-1)$-invertible $\Lambda$-module $M=\Lambda /I$,
where $I$ is generated by $2^{r}$ and $t^n-t+1$, is isomorphic to
$(\mathbb Z/2^r\mathbb Z)^n$ as an abelian group. \qed

\section{$t$-Unipotent $\mathbb Z[t,t^{-1}]$-modules }
\subsection{Properties of $t$-unipotent $\Lambda$-modules}
The following proposition is a simple consequence of Propositions
\ref{cor-krit2} and \ref{cor-krit3}.
\begin{prop} \label{u-cor-krit}
Any $\Lambda$-submodule $M_1$ and any factor module $M/M_1$ of a
Noetherian $(t-1)$-invertible $t$-unipotent $\Lambda$-module $M$
is a $(t-1)$-invertible $t$-unipotent $\Lambda$-module.
\end{prop}

\begin{lem} \label{sum1} Let $M_1,\dots, M_n$ be Noetherian $(t-1)$-invertible
$t$-unipotent $\Lambda$-modules. Then the direct sum
$M=\bigoplus_{i=1}^nM_i$ is a Noetherian $(t-1)$-invertible
$t$-unipotent $\Lambda$-module.
\end{lem}
\proof By Lemma \ref{sum}, $M$ is a Noetherian $(t-1)$-invertible
$\Lambda$-module.

Since $M_i$ is a $(t-1)$-invertible $t$-unipotent
$\Lambda$-module, there is $k_i\in\mathbb N$ such that
$t^{k_i}-1\in \text{Ann}(M_i)$. It is easy to see that $t^k-1\in
\text{Ann}(M)$,  where $k=k_1\dots k_n$, since each polynomial
$t^{k_i}-1$, $i=1,\dots,n$, divides the polynomial $t^k-1$. \qed
\\

Proposition \ref{u-cor-krit} and Lemma \ref{sum1} imply the
following proposition.
\begin{prop} \label{u-cor-krit2}
A Noetherian $\Lambda$-module $M_1$ is $(t-1)$-invertible
$t$-unipotent if and only if each its principle submodule $M_v$ is
$(t-1)$-invertible $t$-unipotent.
\end{prop}

\begin{thm} \label{finit1}
Any Noetherian $\mathbb Z$-torsion $(t-1)$-invertible
$\Lambda$-module is $t$-uni\-potent.
\end{thm}
\proof let $M$ be a Noetherian $\mathbb Z$-torsion
$(t-1)$-invertible $\Lambda$-module. By Corollary \ref{finit}, $M$
consists of finite number of elements. Therefore the automorphism
of $M$, defined by the multiplication by $t$, has a finite order,
say $k$, that is, $t^kv=v$ for all $v\in M$, in other words,
$t^k-1\in \text{Ann}(M)$. \qed \\

The following Propositions \ref{u-biprinciple},
\ref{u-biprinciple1} describe bi-principle $(t-1)$-invertible
$t$-unipotent modules and principle $(t-1)$-invertible
$t$-unipotent modules of mixed type.

\begin{prop} \label{u-biprinciple} Let $M=\Lambda/I$ be a bi-principle
$(t-1)$-invertible $t$-unipotent $\Lambda$-module, and let the
ideal $I=<g(t)>$ is generated by a polynomial $g(t)$. Then
\begin{itemize}
\item[($i$)] all roots of $g(t)$ are roots of unity, \item[($ii$)]
$g(t)$ has not multiple roots, \item[($iii$)] if $\xi$ is a $k$-th
root of unity (that is, $\xi^k=1$), were $k=p^r$ for some prime
$p$, then $\xi$ is not a root of $g(t)$, \item[($iv$)] $g(1)=\pm
1$, \item[($v$)] $\deg g(t)$ is even.
\end{itemize}
\end{prop}
\proof To prove ($i$) and $(ii$), notice that there is $k$ such
that $t^k-1\in I$, since $M$ is $t$-unipotent. Therefore $t^k-1$
is divisible by $g(t)$.

To prove $(iii)$ -- $(v)$, we use Theorem \ref{cor-krit4}. By
Theorem \ref{cor-krit4}, there is a polynomial $f(t)\in I$ such
that $f(1)=1$. We have $f(t)=h(t)g(t)$ for some polynomial
$h(t)\in \mathbb Z[t]$, since $I$ is a principle ideal generated
by $g(t)$. Therefore $g(1)=\pm 1$ (and we can assume that
$g(1)=1$), since we have
$$1=f(1)=h(1)g(1),$$
where $h(1), g(1)\in \mathbb Z$.

On the other hand, if for some prime $p$, a primitive $p^r$-th
root of unity $\xi$ is a root of $g(t)$, then $g(t)$ should be
divided by the $p^r$-th cyclotomic polynomial $\Phi _{p^r}(t)$,
that is, there is a polynomial $h(t)\in \mathbb Z[t]$ such that
$g(t)=\Phi _{p^r}(t)h(t)$. Therefore, $1=g(1)=\Phi _{p^r}(1)h(1)$
and we obtain a contradiction, since $\Phi _{p^r}(1)=p$.

To complete the proof, notice that, by $(iii)$ and $(iv)$,
$\xi=\pm 1$ are not roots of $g(t)$ and hence all roots of $g(t)$
are not real. Thus if $\xi$ is a root of $g(t)$, then the number
$\overline \xi$ complex conjugated to $\xi$  is also a root of
$g(t)$, since $g(t)\in \mathbb Z[t]$. Therefore $\deg g(t)$ is
even, since $\overline \xi\neq \xi$ for all roots of unity $\neq
\pm 1$. \qed

\begin{prop} \label{u-biprinciple1} Let $M=\Lambda/I$ be a principle
$(t-1)$-invertible $t$-unipotent $\Lambda$-module of mixed type,
and let $f(t)=dg(t)$ be the leading generator of the ideal $I$,
where $d\in \mathbb N$ and the polynomial $g(t)$ has content $1$.
Then $g(t)$ satisfies properties $(i)$ -- $(v)$ from Proposition
{\rm \ref{u-biprinciple}}.
\end{prop}
\proof Let $v$ be a generator of $M$. Denote by $M_1$ a
$\Lambda$-submodule of $M$ generated by $v_1=g(t)v$. We have the
exact sequence of $\Lambda$-modules
$$0\to M_1\to M\to M/M_1\to 0,$$
where $M_1$ is a principle module of finite type and $M_2=M/M_1$
is a bi-principle $\Lambda$-module isomorphic to $\Lambda
/<g(t)>$. By Proposition \ref{u-cor-krit}, $M_2$ is
$(t-1)-$invertible $t$-unipotent. Now, we apply Proposition
\ref{u-biprinciple} to complete the proof. \qed \\

Let $M$ be a Noetherian $(t-1)$-invertible $t$-unipotent
$\Lambda$-module. The smallest $k\in \mathbb N$ such that
$t^k-1\in \text{Ann}(M)$ is called the {\it unipotence index} of
$M$.

\begin{lem} \label{divis}
If $M$ is a Noetherian $(t-1)$-invertible $t$-unipotent
$\Lambda$-module of unipotence index $k$, then the polynomial
$\sum_{i=0}^{k-1}t^i\in \text{Ann}(M)$.
\end{lem}
\proof We have $t^k-1=(t-1)(\sum_{i=0}^{k-1}t^i)\in
\text{Ann}(M)$. Therefore $(\sum_{i=0}^{k-1}t^i)v=0$ for all $v\in
M$, since $M$ is a $(t-1)$-invertible $\Lambda$-module. \qed

\begin{lem} \label{dva} A Noetherian $(t-1)$-invertible
$\Lambda$-module $M$ of unipotence index $2$ is a finite $\mathbb
Z$-module of odd order.
\end{lem}
\proof It follows from Lemma \ref{p^r} and Corollary \ref{finit}
that $M$ is finite. By Lemma \ref{divis}, the polynomial $(t+1)\in
\text{Ann}(M)$. Therefore $tv=-v$ for all $v\in M$. In particular,
if $v$ is of order $2$, then $tv=v$. This is impossible, since $M$
is $(t-1)$-invertible. Therefore $M$ has not elements of even
order. \qed
\begin{prop} \label{cyclicgroup} A cyclic group $G$ of order $n=p_1^{r_1}\dots
p_m^{r_m}$, where $p_1,\dots,p_m$ are primes, possesses a
structure of $(t-1)$-invertible $\Lambda$-module of unipotence
index $k$ if and only if for each  $i=1,\dots,m$ the polynomial
$\sum_{i=0}^{k-1}t^i$ has a root $a_i\neq 1$ in the field $\mathbb
Z/p_i\mathbb Z$.
\end{prop}
\proof By Theorem \ref{decomposition}, it suffices to consider
only the case when $i=1$, that is $n=p^m$ for some prime $p$.

Let a cyclic group $G$ of order $n=p^m$ has a structure of
$(t-1)$-invertible $\Lambda$-module of unipotence index $k$, then
its subgroup $G_p=p^{m-1}G$ consisting of the elements of order
$p$ is also a $(t-1)$-invertible $\Lambda$-module of unipotence
index $k$. Therefore the polynomial $\sum_{i=0}^{k-1}t^i\in
\text{Ann}(G_p)$. Let $v\in G_p$ be a generator of $G_p$, then
$tv=av$ for some $a \not\equiv 1\mod p$ since $G_p$ is a
$(t-1)$-invertible module. We have $\sum_{i=0}^{k-1}a^iv=0$.
Therefore $\sum_{i=0}^{k-1}a^i\equiv 0\mod p$, that is, the
polynomial $\sum_{i=0}^{k-1}t^i$  has a root in the field $\mathbb
Z/p\mathbb Z$ not equal to $0$ or $1$.

Conversely, let $a\not\equiv 1\mod p$ be a root of the polynomial
$\sum_{i=0}^{k-1}t^i$ in the field $\mathbb Z/p_i\mathbb Z$, and
let $v$ be a generator of a cyclic group $G$ of order $p^r$. If we
define the action of $t$ on the $\mathbb Z$-module $G$ putting
$t(v)=av$, we obtain a structure of $(t-1)$-invertible
$\Lambda$-module on $G$, since $a\not\equiv 1\mod p$. It is easy
to see that $t^k-1\in \text{Ann}(G)$. \qed

\begin{thm} \label{zfinite}
Any Noetherian  $(t-1)$-invertible $t$-unipotent $\Lambda$-module
$M$ is finitely generated over $\mathbb Z$.
\end{thm}
\proof Theorem follows from Proposition \ref{finitenes}, since for
some $k\in \mathbb Z$ the polynomial $t^k-1\in \text{Ann}(M)$. \qed \\

It follows from Theorem \ref{finit1} and Structure Theorem for
finitely generated $\mathbb Z$-modules that a Noetherian
$(t-1)$-invertible $t$-unipotent $\Lambda$-module $M$ as a
$\mathbb Z$-module is isomorphic to
\begin{equation}
\label{decomp} M\simeq M_{fin}\oplus \mathbb Z^k,
\end{equation}
where $M_{fin}$ is the submodule of $M$ consisting of the elements
of finite order. The rank $k$ of the free part of $M$ in
decomposition (\ref{decomp}) is called {\it Betti number} of
Noetherian $(t-1)$-invertible $t$-unipotent $\Lambda$-module $M$.
\begin{thm} \label{even} The Betti number of a Noetherian
$(t-1)$-invertible $t$-unipotent $\Lambda$-module $M$ is an even
number.
\end{thm}
\proof By definition, the Betti number of $M$ coincides with Betti
number of the Noetherian $(t-1)$-invertible $t$-uni\-po\-tent
$\Lambda$-module $M_{free}=M/M_{fin}$.

The module $M_{free}$ has not non-zero elements of finite order.
Therefore the annihilator $\text{Ann}_v$ of each its element $v$
is a principle ideal, it is generated by polynomial $g_v(t)$
satisfying properties $(i)$ -- $(v)$  from Proposition
\ref{u-biprinciple}.

Let $M_{free}$ is generated by  elements $v_1,\dots,v_m$ over
$\Lambda$. Then there is a surjective $\Lambda$-homomorphism
$$f:\Lambda/<g_{v_1}(t)>\oplus \dots \oplus \Lambda/<g_{v_m}(t)>\longrightarrow M_{free}.$$

Consider the modules $\widetilde M=\bigoplus \Lambda/<g_{v_i}(t)>$
and $M_{free}$ as free $\mathbb Z$-modules and denote by
$h_{\widetilde M}$ and $h_{M_{free}}$ the automorphisms
respectively of $\widetilde M$ and $M_{free}$ defined by the
multiplication by $t$. Then it is easy to see that the
characteristic polynomial $\widetilde \Delta(t)=\det
(h_{\widetilde M}-t\text{Id})$ coincides up to the sign with the
product $g_{v_1}(t)\dots g_{v_m}(t)$. Next,  the characteristic
polynomial $\Delta(t)=\det (h_{ M_{free}}-t\text{Id})$ is a
divisor of the polynomial $\widetilde \Delta(t)$, since the
homomorphism $f$ is surjective and $t$-equivariant. Therefore all
roots of $\Delta(t)$ are roots of unity $\neq \pm 1$ and hence
$\deg \Delta(t)$ is an even number. To complete the proof, notice
that the Betti number of $M_{free}$ coincides with $\deg
\Delta(t)$. \qed
\subsection{Derived Alexander modules}
To a Noetherian $(t-1)$-invertible $\Lambda$-module $M$ we
associate an infinite sequence of  Noetherian $(t-1)$-invertible
$t$-unipotent $\Lambda$-modules
\begin{equation} \label{seque}
A_n(M)=M/(t^n-1)M, \qquad n\in \mathbb N.
\end{equation}
The module $A_n(M)$ is called the {\it $n$-th derived Alexander
module} of $\Lambda$-module $M$.

Note that $A_1(M)=0$, since $M$ is $(t-1)$-invertible. It is also
evident that $A_n(A_n(M))=A_n(M)$.

It is obvious, that if $f: M_1\to M_2$ is a $\Lambda$-homomorphism
of $(t-1)$-invertible modules, then the sequence of
$\Lambda$-homomorphisms $$f_{n*}:A_n(M_1)\to A_n(M_2),$$ $n\in
\mathbb N$, is well defined, that is, the map $M\mapsto \{
A_n(M)\}$ is a functor from the category of Noetherian
$(t-1)$-invertible $\Lambda$-modules to the category of infinite
sequences of Noetherian $(t-1)$-invertible $t$-unipotent
$\Lambda$-modules.
\begin{prop} \label{simple} If
$$0\to M_1\stackrel{f}\longrightarrow M\stackrel{g}\longrightarrow M_2\to 0$$
is an exact sequence of Noetherian $(t-1)$-invertible
$\Lambda$-modules, then $$A_n(M_2)\simeq A_n(M)/\text{im}\,
f_{n*}(A_n(M_1)).$$

If $M=\bigoplus_{i=1}^kM_i$ is the direct sum of  Noetherian
$(t-1)$-invertible $\Lambda$-modules $M_i$, then $$A_n(M)\simeq
\bigoplus_{i=1}^kA_n(M_i).$$
\end{prop}
\proof Obvious. \qed
\begin{prop} \label{p^r-1} Let $p$ be a prime number and $r\in \mathbb N$,
then for a Noetherian $(t-1)$-invertible $\Lambda$-module $M$ its
derived Alexander module $A_{p^r}(M)$ is finite.
\end{prop}
\proof It follows from Lemma \ref{p^r} and Corollary \ref{finit}.
\qed

\begin{ex} \label{ex1} For $M_m=\Lambda/\langle (m+1)t-m\rangle$,
where $m\in\mathbb N$, its $n$-th derived Alexander module
$$A_n(M_m)\simeq \mathbb Z/((m+1)^n-m^n)\mathbb
Z$$ is a cyclic group of order $(m+1)^n-m^n$ and the
multiplication by $t$ is given by
$$tv=(-1)^{n+1}m(\sum_{i=1}^{n-1}(-1)^{i}{n \choose i}(m+1)^{n-i-1})v$$ for
all $v\in A_n(M_m)$.
\end{ex}
\proof The module $M_m=\Lambda/\langle (m+1)t-m\rangle$ is
isomorphic to a $\Lambda$-submodule $\mathbb
Z[\frac{m}{m+1},\frac{m+1}{m}] \subset \mathbb Q$ if we put
$t=\frac{m}{m+1}$ and $tv=\frac{m}{m+1}v$ for $v\in \mathbb Q$.
Therefore we have
\begin{center} $A_n(M_m)\simeq M_m/(t^n-1)M_m\simeq \mathbb
Z[\frac{m+1}{m},\frac{m}{m+1}]/\langle(\frac{m}{m+1})^n-1\rangle$
\end{center}
and consequently,
\begin{center} $A_n(M_m)\simeq \mathbb Z[\frac{m}{m+1},\frac{m+1}{m}]/\langle(m+1)^n-m^n\rangle .$
\end{center}

It is easy to see that the module $\mathbb
Z[\frac{m}{m+1},\frac{m+1}{m}]$ coincides with the sum of
submodules $\mathbb Z[\frac{1}{m+1}]$ and $\mathbb
Z[\frac{1}{m}]\subset \mathbb Q$,
\begin{center} $\mathbb Z[\frac{m}{m+1},\frac{m+1}{m}]=\mathbb
Z[\frac{1}{m+1}]+\mathbb Z[\frac{1}{m}].$
\end{center}
Indeed, it is obvious, that
\begin{center} $\mathbb Z[\frac{m}{m+1},\frac{m+1}{m}]\subset \mathbb
Z[\frac{1}{m+1}]+\mathbb Z[\frac{1}{m}].$
\end{center}
Next, we have
\begin{center}$(\frac{m+1}{m})^n=\frac{\sum_{i=0}^n{n \choose i}
m^{n-i}}{m^n}$ \end{center} and therefore
\begin{center}
$\frac{1}{m^n}=(\frac{m+1}{m})^n -\sum_{i=0}^{n-1}{n \choose
i}\frac{1}{m^i}$.\end{center} Similarly, we have
\begin{center}
$\frac{1}{(m+1)^n}=\sum_{i=0}^{n-1}(-1)^{n+1+i}{n \choose
i}\frac{1}{(m+1)^i}+(-1)^n(\frac{m}{m+1})^n $ \end{center} In
particular, $\frac{1}{m}=\frac{m+1}{m}-1$ and
$\frac{1}{m+1}=1-\frac{m}{m+1}$. Therefore, by induction, we
obtain that $\frac{1}{m^n}$, $\frac{1}{(m+1)^n}\in \mathbb
Z[\frac{m}{m+1},\frac{m+1}{m}]$ for all $n$ and hence
\begin{center} $\mathbb
Z[\frac{1}{m+1}]+\mathbb Z[\frac{1}{m}]\subset \mathbb
Z[\frac{m}{m+1},\frac{m+1}{m}].$
\end{center}

Let us show now that each element $v\in \mathbb
Z[\frac{m}{m+1},\frac{m+1}{m}]$ is equivalent to some $v_{in}\in
\mathbb Z\subset \mathbb Z[\frac{m}{m+1},\frac{m+1}{m}]$ modulo
the ideal $I=\langle(m+1)^n-m^n\rangle$. For this, it suffices to
show that for each $k$ there are $v_{in,k},u_{in,k}\in \mathbb Z$
such that \begin{center} $\frac{1}{m^k}\equiv v_{in,k}\mod I\qquad
\text{and}\qquad \frac{1}{(m+1)^k}\equiv u_{in,k}\mod I$.
\end{center}

We prove the existence of such elements only for $\frac{1}{m^k}$
and the case $\frac{1}{(m+1)^k}$ will be left to the reader, since
it is similar. We have
\begin{center} $\frac{(m+1)^n-m^n}{m^{k}}=\sum_{i=1}^{n} {n
\choose i}m^{n-i-k}\equiv 0\mod I$
\end{center} and therefore
\begin{center} $\frac{1}{m^{k}}\equiv -\sum_{j=k+1-n}^{k-1} {n \choose
n+j-k}\frac{1}{m^{j}}\mod I$.
\end{center}
In particular,
\begin{center} $\frac{1}{m}\equiv -\sum_{j=0}^{n-2} {n \choose
n-j-1}m^{j}\mod I$.
\end{center}
Now the existence of desired $v_{in,k}$ is proved by induction on
$k$.

It follows from the above consideration that
\begin{center} $A_n(M_m)\simeq \mathbb Z[\frac{m}{m+1},\frac{m+1}{m}]/\langle(m+1)^n-m^n\rangle $
\end{center}
is a cyclic group generated by the image $\overline 1$ of $1\in
\mathbb Z[\frac{m}{m+1},\frac{m+1}{m}]$. We have
$$((m+1)^n-m^n)\overline 1=0$$ and hence the order of $A_n(M_m)$ is
a divisor of $(m+1)^n-m^n$.

Let us show that the order of $A_n(M_m)$ is equal to
$(m+1)^n-m^n$. Let $k\in \mathbb Z$ be such that $k\overline 1=0$.
Then \begin{center} $k=(\sum_{i_1\leq i\leq
i_2}a_i\frac{1}{(m+1)^i}+ \sum_{j_1\leq j\leq
j_2}b_j\frac{1}{m^j})((m+1)^n-m^n)$,
\end{center}
where $a_i, b_j\in\mathbb Z$. Multiplying by $(m+1)^{i_2}$ and
$m^{j_2}$ if $i_2>0$ or $j_2>0$, we obtain an equality
$$(m+1)^{i_2}m^{j_2}k =C((m+1)^n-m^n)$$
with some $C\in\mathbb Z$ which shows that $(m+1)^n-m^n$ is a
divisor of $k$, since $m$, $m+1$, and $(m+1)^n-m^n$ are coprime.

To calculate the action of $t$ on the cyclic group
$$A_n(M_m)\simeq \mathbb Z/((m+1)^n-m^n)\mathbb Z,$$ notice that
\begin{center}
$t\overline 1=\overline{
\frac{m}{m+1}}=(-1)^{n+1}m(\sum_{i=1}^{n-1}(-1)^{i}{n \choose
i}(m+1)^{n-i-1})\overline 1,$
\end{center}
since similar (as above) calculation gives
\begin{center}
$\qquad \frac{1}{m+1}\equiv (-1)^{n+1}\sum_{i=1}^{n-1}(-1)^{i}{n
\choose i}(m+1)^{n-i-1} \mod I.\qquad \qquad \qed $
\end{center}

\begin{prop} \label{realization}
An abelian group $G$ is isomorphic {\rm (}as a $\mathbb
Z$-module{\rm )} to the derived Alexander module $A_2(M)$ of some
Noetherian $(t-1)$-invertible $\Lambda$-module $M$ if and only if
$G$ is a finite group of odd order.
\end{prop}
\proof By Lemma \ref{dva}, we need only to prove that for any
finite group $G$ of odd order there is a  Noetherian
$(t-1)$-invertible $\Lambda$-module $M$ for which $A_2(M)\simeq
G$.

Represent $G$ as a direct sum of cyclic groups:
$$G=\bigoplus_{i=1}^k G_i,$$
and let $n_i=2m_i+1$ be the order of $G_i$.

For each $i$, consider the $\Lambda$-module $M_{m_i}$ from Example
\ref{ex1}. We have $A_2(M_{m_i})$ is a cyclic group of order
$(m_i+1)^2-m_i^2=2m_i+1$. Now, proposition follows from
Proposition \ref{simple} if we put $M=\bigoplus_{i=1}^kM_{m_i}$.
\qed

\begin{thm} \label{period} Let $M$ be a Noetherian $(t-1)$-invertible
$t$-unipotent $\Lambda$-module of unipotence index $k$. Then the
sequence of its derived Alexander modules $$A_1(M), \dots, A_n(M),
\dots $$ has period $k$, that is, $A_n(M)\simeq A_{n+k}(M)$ for
all $n$.

If $n$ and $k$ are coprime, then $A_n(M)=0$.
\end{thm}
\proof Note that if $k$ is the unipotence index of $M$, then, by
Lemma \ref{divis}, the polynomial $f_k(t)=\displaystyle
\sum_{i=0}^{k-1}t^i\in \text{Ann}(M)$. Besides, to get $A_n(M)$
from $M$, it suffices to factorize $M$ by the relations
$f_n(t)v=0$ for all $v\in M$, where $f_n(t)=\displaystyle
\sum_{i=0}^{n-1}t^i$. Now, to prove the periodicity of sequence
(\ref{seque}), it suffices to notice that
$$f_{n+k}(t)=t^nf_k(t)+f_n(t).$$

Let $n$ and $k$ be coprime and let polynomials $f_k(t)$ and
$f_n(t)$ belong to $\text{Ann}(M)$. Applying Euclidian algorithm
to $f_k(t)$ and $f_n(t)$, it is easy to see that there are
polynomials $g_k(t)$ and $g_n(t)$ such that
$$f_k(t)g_k(t)+f_n(t)g_n(t)=1,$$
since $n$ and $k$ are coprime. Therefore $\text{Ann}(M)=\Lambda$
and hence $A_n(M)=0$. \qed
\begin{ex} \label{ex2}
The $\Lambda$-module $M=\Lambda/<t^2-t+1>$ has the following
derived Alexander modules:
$$A_{6k\pm 1}(M)=0,\qquad A_{6k\pm 2}(M)\simeq \mathbb Z/3\mathbb
Z,\qquad A_{6k+3}(M)\simeq (\mathbb Z/2\mathbb Z)^2, $$ where the
multiplication by $t$ on $\mathbb Z/3\mathbb Z$ coincides with the
multiplication by $2$ and  the multiplication by $t$ on $(\mathbb
Z/2\mathbb Z)^2$ coincides with cyclic permutation of the non-zero
elements of $A_{6k+3}(M)$.
\end{ex}
\proof The module $M$ has the unipotency index $6$, since
$t^2-t+1$ is a divisor of the polynomial $t^6-1$. Therefore
$A_{6k\pm 1}(M)=0$.

To compute $A_{6k+2}(M)$, it suffices to compute $A_2(M)$. We have
$A_2(M)=\Lambda/<t^2-t+1, t+1>$ and since
$$t^2-t+1=(t-2)(t+1)+3,$$
then $\Lambda/<t^2-t+1, t+1>=\Lambda/<t+1, 3>\simeq \mathbb
Z/3\mathbb Z$.

To compute $A_{6k+3}(M)$, it suffices to compute $A_3(M)$. We have
$A_3(M)=\Lambda/<t^2-t+1, t^2+t+1>$ and since
$$t^2+t+1=t^2-t+1+2t,$$
then $\Lambda/<t^2-t+1, t^2+t+1>=\Lambda/<t^2+t+1, 2>\simeq
\mathbb (Z/2\mathbb Z)^2$.

To compute $A_{6k+4}(M)$, it suffices to compute $A_4(M)$. We have
$A_4(M)=\Lambda/<t^2-t+1, t^3+t^2+t+1>$ and since
$$t^3+t^2+t+1=(t+2)(t^2-t+1)+2t-1,$$
then $\Lambda/<t^2-t+1, t^3+t^2+t+1>=\Lambda/<t^2-t+1, 2t-1>$ is
isomorphic to the quotient module $M/(2t-1)M$. Let $v$ be a
generator of bi-principle module $M$. It is easy to check that in
the basis $v_1=v,v_2=tv$ of $M$ over $\mathbb Z$, the module
$(2t-1)M$ is generated by the elements $2v_2-v_1$ and
$t(2v_2-v_1)=v_2-2v_1$, since $tv_2=v_2-v_1$. In the new basis
$e_1=v_1$, $e_2=v_2-2v_1$, the element $2v_2-v_1=2e_2+3e_1$, that
is, $(2t-1)M$ is generated over $\mathbb Z$ by $3e_1$ and $e_2$.
Therefore $A_4(M)\simeq \mathbb Z/3\mathbb Z$. \qed

\section{Alexander modules of irreducible $C$-groups}
\subsection{Proof of Theorems \ref{main1} and \ref{main2}}
\label{sec3} Recall that the class of irreducible $C$-groups
coincides with the class of fundamental groups of knotted
$n$-manifolds $V$ if $n\geq 2$ and the knot groups are also
$C$-groups if they are gvien by Wirtinger presentation. Similarly,
the class of irreducible Hurwitz $C$-groups coincides with the
class of the fundamental groups of the complements of
irreducible\, "affine" Hurwitz (resp., pseudo-holomorphic) curves
and it contains the subclass of the fundamental groups of the
complements of algebraic irreducible affine plane curves.
Therefore to speak about the Alexander modules of knotted
$n$-manifolds and, respectively, about the Alexander modules of
irreducible Hurwitz (resp., pseudo-holomorphic) curves is the same
as to speak about the Alexander modules of irreducible $C$-groups
and, respectively, of irreducible Hurwitz $C$-groups. Hence
Theorems \ref{main1} and \ref{main2} are equivalent to the
following two theorems.

\begin{thm} \label{main11} A  $\Lambda$-module $M$ is the Alexander
module of an irreducible $C$-group if and only if it is
Noetherian $(t-1)$-invertible.
\end{thm}

\begin{thm} \label{main22} A  $\Lambda$-module $M$ is the Alexander
module of an irreducible Hurwitz $C$-group if and only if it is
Noetherian $(t-1)$-invertible $t$-unipotent $\Lambda$-module.

The unipotence index of the Alexander module $A_0(G)$ of an
irreducible $C$-group $G$ of degree $m$ is a divisor of $m$.
\end{thm}

\proof Let
\begin{equation} \label{pres} G=<x_1,\dots,x_m\, \mid\,
r_1,\dots,r_n>\end{equation} be a $C$-presentation of a $C$-group
$G$ and $\mathbb F_m$ be the free group freely generated by the
$C$-generators $x_1,\dots,x_m$. Denote by
$\frac{\partial}{\partial x_i}$ the Fox derivative (\cite{C-F}),
that is, an endomorphism of the group ring $\mathbb Z[\mathbb
F_m]$ over $\mathbb Z$ of the free group $\mathbb F_m$ into
itself, such that $\frac{\partial}{\partial x_i} :\mathbb
Z[\mathbb F_m]\to \mathbb Z[\mathbb F_m]$ is a $\mathbb Z$-linear
map defined by the following properties
\begin{equation} \label{fox}
 \begin{array}{rcl}
\displaystyle \frac{\partial x_j}{\partial x_i} & = & \delta_{i,j} \\
\displaystyle \frac{\partial uv}{\partial x_i} & = & \displaystyle
\frac{\partial u}{\partial x_i}+u\frac{\partial v}{\partial x_i}
\end{array}
\end{equation}
for any $u,v\in \mathbb Z[\mathbb F_m]$.  The matrix
$$\mathcal A(G)=\nu_*\bigl(\tfrac{\partial r_i}{\partial x_j}\bigr)\in \text{Mat}_{n\times
m}(\mathbb Z[t,t^{-1}])$$ is called the {\it Alexander matrix} of
the $C$-group $G$ given by presentation (\ref{pres}), where $r_i$,
$i=1,\dots,n$, are the defining relations of $G$ and $\nu_*:
\mathbb Z[\mathbb F_m]\to \mathbb Z[\mathbb F_1]\simeq \mathbb
Z[t,t^{-1}]$ is induced by the canonical $C$-epimorphism $\nu
:\mathbb F_m\to \mathbb F_1$.

\begin{lem} \label{zero} The sum of the columns of the Alexander
matrix $\mathcal A(G)$ of a $C$-group $G$, given by presentation
{\rm (\ref{pres}\rm )}, is equal to zero.
\end{lem}
\proof Each relation $r_i$ has the form
$$r= wx_jw^{-1}x_l^{-1},$$
where $w$ is a word in letters $x_1^{\pm 1},\dots,x_m^{\pm 1}$,
and $x_j,x_l$ are some two letters.

By induction on the length $l(w)$ of the word $w$, let us show
that $$\sum_{k=1}^m\nu_*(\frac{\partial r}{\partial x_k})=0.$$ If
$l(w)=0$, that is, $r:=x_jx_l^{-1}$, we have
$$\nu_*(\frac{\partial r}{\partial x_k})=
\left\{
\begin{array}{cl}
1\qquad & \text{if}\, \, k=j, \\
-1\qquad & \text{if}\, \, k=l, \\
0\qquad & \text{if}\, \, k\neq j\, \,\text{and}\,\, k\neq l
\end{array}
\right.
$$
and in this case we obtain $\sum_{k=1}^m\nu_*(\frac{\partial
r}{\partial x_k})=0$.

Assume that for all words $r= wx_jw^{-1}x_l^{-1}$ we have
$\sum_{k=1}^m\nu_*(\frac{\partial r}{\partial x_k})=0$ if
$l(w)\leq n$. Consider a word $r= wx_jw^{-1}x_l^{-1}$, such that
$l(w)=n+1$. Put $r_1= w_1x_jw_1^{-1}x_l^{-1}$, where
$w=x_i^{\varepsilon}w_1$, $\varepsilon=\pm 1$, and $l(w_1)=n$. We
consider only the case when $i\neq j$, $i\neq l$, $j\neq l$, and
$\varepsilon =1$. All other cases are similar and the proof that
$\sum_{k=1}^m\nu_*(\frac{\partial r}{\partial x_k})=0$ in these
cases will be left to the reader.

It follows from (\ref{fox}) that
$$\nu_*(\frac{\partial r}{\partial x_k})=
\left\{
\begin{array}{ll} \displaystyle
t\nu_*(\frac{\partial r_1}{\partial x_k})\qquad & \text{if}\, \, k\neq i,k\neq j,k\neq l, \\
 \displaystyle 1+t\nu_*(\frac{\partial r_1}{\partial x_k})-t\qquad & \text{if}\, \, k=i, \\
 \displaystyle t\nu_*(\frac{\partial r_1}{\partial x_k})\qquad & \text{if}\, \, k=j, \\
 \displaystyle t(\nu_*(\frac{\partial r_1}{\partial x_k})+1)-1\qquad & \text{if}\, \, k=l
\end{array}
\right.
$$
and it is easy to see that $\sum_{k=1}^m\nu_*(\frac{\partial
r}{\partial x_k})=0$. \qed \\

To each monomial $a_it^i\in \mathbb Z[t]$ let us associate a word
$$w_{a_it^i}(x_1,x_2)=(x_{2}^ix_1x_{2}^{-(i+1)})^{a_i}$$ if
$a_i>0$ and
$$w_{a_it^i}(x_1,x_{2})=(x_{2}^{i+1}x_1^{-1}x_{2}^{-i})^{-a_i}$$
if $a_i<0$, and for $g(t)=\sum_{i=0}^k a_it^i\in \mathbb Z$ we put
$$w_{g(t)}(x_1,x_{2})=\prod_{i=0}^k w_{a_it^i}(x_1,x_{2}).$$
Next, to a polynomial $f(t)=(1-t)g(t) +1$ we associate a word
\begin{equation} \label{rel10}
r_{f(t)}(x_1,x_2)=w_{g(t)}(x_1,x_{2})x_1w^{-1}_{g(t)}(x_1,x_2)x_2^{-1},
\end{equation}
and to a vector $u=(1-t)\overline u=((1-t)g_{1}(t),\dots
,(1-t)g_{m}(t))$, we associate a word
\begin{equation} \label{rel11}
r_u(x_1,\dots,x_{m+1})=w_u(x_1,\dots,x_{m+1})x_{m+1}w^{-1}_u(x_1,\dots,x_{m+1})x^{-1}_{m+1},
\end{equation}
where
$$w_u(x_1,\dots, x_{m+1})=\prod_{i=1}^mw_{g_i(t)}(x_i,x_{m+1}).$$

\begin{lem} \label{ftor} For a polynomial $f(t)=(1-t)g(t)+1$ and a vector
$$u=((1-t)g_{1}(t),\dots ,(1-t)g_{m}(t))$$ we have
$$\begin{array}{l}
\displaystyle \nu_*(\frac{\partial r_{f(t)}}{\partial x_1})=f(t),
\\ \displaystyle \nu_*(\frac{\partial r_{u}}{\partial
x_i})=(1-t)g_i(t), \qquad i=1,\dots,m.
\end{array}
$$
\end{lem}
\proof Let $f(t)=(1-t)g(t)+1$. It follows from (\ref{fox}) that
$$\nu_*(\frac{\partial w_{g(t)}(x_1,x_{2})}{\partial
x_1})=-\nu_*(\frac{\partial w^{-1}_{g(t)}(x_1,x_{2})}{\partial
x_1})=g(t)$$ since
$w_{g(t)}(x_1,x_{2})w^{-1}_{g(t)}(x_1,x_{2})=1$,
$$\nu_*(w_{g(t)}(x_1,x_{2}))=\nu_*(w_{a_it^i}(x_1,x_{2}))=1,$$ and
$$\nu_*(\frac{\partial w_{a_it^i}(x_1,x_{2})}{\partial
x_1})=a_it^i.$$ Therefore we have
$$\begin{array}{ll}
\displaystyle \nu_*(\frac{\partial r_{f(t)}}{\partial x_1})= &
\displaystyle \nu_*(\frac{\partial
(w_{g(t)}(x_1,x_{2})x_1w^{-1}_{g(t)}(x_1,x_{2})x_{2}^{-1})}{\partial
x_1})= \\ & g(t)+1-tg(t)=f(t).
\end{array}$$

The proof in the second case  is similar  and it will be left to
the reader. \qed

\begin{prop} \label{kum}
The Alexander module $A_0(G)$ of a $C$-group $G$, given by
presentation {\rm (\ref{pres})}, is isomorphic to a factor module
$\Lambda^{m-1}/M(G)$, where the submodule $M(G)$ of
$\Lambda^{m-1}$ is generated by the rows of the matrix
$\overline{\mathcal A}$ formed by the first $m-1$ columns of the
Alexander matrix $\mathcal A(G)$.
\end{prop}
\proof To describe the Alexander module of a $C$-group $G$, we
follow \cite{M} (see also \cite{Ku3}). To a $C$-group $G$ given by
$C$-presentation (\ref{pres}) we associate a two-dimensional
complex $K$ with a single vertex $x_0$ whose one dimensional
skeleton is a bouquet of oriented circles $s_i$, $1\leq i\leq m$,
corresponding to the $C$-generators of $G$ in presentation
(\ref{pres}). Furthermore, $K\setminus (\cup
s_i)=\bigsqcup_{i=1}^n\stackrel{\circ}D_i$ is a disjoint union of
open discs. Each disc $D_i$ corresponds to the relation
$r_i=x_{j_{i,1}}^{\varepsilon_{i,1}}\dots
x_{j_{i,k_i}}^{\varepsilon_{i,k_i}}$ from presentation
(\ref{pres}), where $\varepsilon_{i,j}=\pm 1$, and it is glued to
the bouquet $\bigvee s_i$ along the path
$s_{j_{i,1}}^{\varepsilon_{i,1}}\dots
s_{j_{i,k_i}}^{\varepsilon_{i,k_i}}$. It is clear that
$\pi_1(K,x_0)\simeq G$.

The $C$-homomorphism $\nu :G\to \mathbb F_1$ defines an infinite
cyclic covering $f:\widetilde K\to K$ such that $\pi_1(\widetilde
K)=N$ and $H_1(\widetilde K,\mathbb Z)=N/N'$, where $N=\ker \nu$.
The group $\mathbb F_1$ acts on $\widetilde K$.

Let $\widetilde K_0=f^{-1}(x_0)$, and let $\widetilde K_1$ be the
one-dimensional skeleton of the complex $\widetilde K$. Consider
the following exact sequences of homomorphisms of homology groups
with coefficients in $\mathbb Z$:
\begin{equation} \label{exact}
\begin{CD}
@. @. @. 0 @.
\\
@. @.@. @VVV @.
\\
@.
@. @. 
H_1(\widetilde K) 
@.
\\
@. @.@. @VVV @.
\\
@>>> H_2(\widetilde K,\widetilde K_1) @>{\alpha}>> H_1(\widetilde
K_1,\widetilde K_0)  @>{\beta}>>
 H_1(\widetilde K,\widetilde K_0) @>>> 0
  \\
 @.
 @.@. @VV{\partial}V @.
 \\
 @.
 @. @. H_0(\widetilde K_0) @.
 \\
@. @.@. @VVV @.
\\
@. @. @. 0 @.
\end{CD}
\end{equation}
The action of $\mathbb F_1$ on $\widetilde K$ turns the groups in
these sequences into $\Lambda$-modules. We fix a vertex $p_0\in
\widetilde K_0$. Let $p_i=t^ip_0$ be the result of action of the
element $t^i\in \mathbb F_1$ on the point $p_0$. Then
$H_1(\widetilde K_1,\widetilde K_0)$ is a free $\Lambda$-module
whose generators $\bar s_i$ are edges joining $p_0$ with $p_1$
which are mapped onto the loops $s_i$. The result of action of
$t^i$ on the generator $\bar s_j$ is an edge beginning at the
vertex $p_i$ which is mapped onto the loop $s_j$.

The free $\Lambda$-module $H_2(\widetilde K,\widetilde K_1)$ is
generated by the discs $\overline D_i$, $i=1,\dots,n$,
corresponding to the relations
$r_i=x_{j_{i,1}}^{\varepsilon_{i,1}}\dots
x_{j_{i,k_i}}^{\varepsilon_{i,k_i}}$, where each disc $\overline
D_i$ is glued to the one-dimensional skeleton along the product of
paths  $$\displaystyle t^{\delta(\varepsilon_{i,1})}\bar
s_{j_{i,1}}^{\varepsilon_{i,1}},
t^{\delta(\varepsilon_{i,2})+\varepsilon_{i,1}}\bar
s_{j_{i,2}}^{\varepsilon_{i,2}}, \dots,
t^{\delta(\varepsilon_{i,k_i})+\sum_{l=1}^{k_i-1}\varepsilon_{i,l}}\bar
s_{j_{i,k_i}}^{\varepsilon_{i,k_i}},$$ where $\delta(1)=0$ and
$\delta(-1)=-1$. It is easy to verify that the coordinates of
elements $\alpha(\overline D_i)\in H_1(\widetilde K_1,\widetilde
K_0)$ in the basis $\bar s_1,\dots,\bar s_m$ coincide with the
rows $\mathcal A_i$ of the Alexander matrix $\mathcal A(G)$ of
$C$-group $G$ given by presentation (\ref{pres}).

It follows from the vertical exact sequence in (\ref{exact}) that
$\partial(\beta(\bar s_i))=(t-1)p_0$ for each generator $\bar s_i$
of the module $H_1(\widetilde K_1,\widetilde K_0)$. Let us choose
a new basis $e_i=\bar s_i-\bar s_m$, $i=1,\dots,m-1$, $e_m=\bar
s_m$ in $H_1(\widetilde K_1,\widetilde K_0)$. Then $\beta(e_i)\in
\ker
\partial$ for $i=1,\dots,m-1$, and $\ker \partial$ is generated by
$\beta(e_1),\dots, \beta(e_{m-1})$. Hence we may identify
$H_1(\widetilde K)$ with $\beta(H'_1(\widetilde K_1,\widetilde
K_0))$, where $H'_1(\widetilde K_1,\widetilde K_0)$ is a free
submodule of the free $\Lambda$-module $H_1(\widetilde
K_1,\widetilde K_0)$ generated by the elements
$e_1,\dots,e_{m-1}$.

In the basis $e_1,\dots,e_m$ the matrix formed by the coordinates
of $\alpha(\overline D_i)$ coincides with the matrix
$\widetilde{\mathcal A}(G)$ obtained from $\mathcal A(G)$ by
replacing the last column by the column of zeros. Hence
$H_1(\widetilde K)$ is isomorphic to the quotient of the free
$\Lambda$-module $H'_1(\widetilde K_1,\widetilde K_0)\simeq
\bigoplus_{i=1}^{m-1}\Lambda e_i$ by the submodule $M(G)$
generated by the rows of the matrix $\overline{\mathcal A}(G)$,
where $\overline{\mathcal A}(G)$ is the matrix formed by the first
$m-1$ columns of the matrix $\mathcal A(G)$. \qed \\

To prove that a  Noetherian $(t-1)$-invertible (resp.,
$t$-unipotent) $\Lambda$-module $M$ is the Alexander module of an
irreducible (resp., Hurwitz) $C$-group, we use Proposition
\ref{usef}. By Proposition \ref{usef}, a Noetherian
$(t-1)$-invertible $\Lambda$-module $M$ is isomorphic to a factor
module $\Lambda^m/M_1$ of a free $\Lambda$-module $\Lambda^m$,
where the submodule $M_1$ is generated by elements
$u_1,\dots,u_m,\dots,u_{m+k}$ of $\Lambda^m$ such that
\begin{itemize}
\item[($i$)] for  $i=1,\dots, m$ the vector $u_i= (0, \dots, 0,
f_i(t),0,\dots,0)$, where a polynomial $f_i(t)$  is such that
$f_i(1)=1$ and it stands on the $i$-th place, \item[($ii$)]
$u_{m+j}=(1-t)\overline u_{m+j}=((1-t)g_{j,1}(t),\dots
,(1-t)g_{j,m}(t))$ for $j=1,\dots,k$, where $g_{j,l}(t)$ are
polynomials,
\end{itemize}
and if $M$ is a $t$-unipotent $\Lambda$-module of unipotence index
$n$, then we can assume that
\begin{itemize}
\item[($iii$)] the vector $u_{m+k+i}= (0, \dots,
0,t^n-1,0,\dots,0)\in M_1$  for $i=1,\dots, m$, where the
polynomial $t^n-1$ stands on the $i$-th place.
\end{itemize}

Express each polynomial $f_i(t)$ in the form
$f_i(t)=(1-t)g_i(t)+1$ and consider a $C$-group
$$G=\langle x_1,\dots,x_{m+1} \mid r_1,\dots,r_{m+k}\rangle ,$$
where $r_i:=r_{f_i(t)}(x_i,x_{m+1})$ for $i=1,\dots,m$ and
$r_{m+j}:=r_u(x_1,\dots, x_{m+1})$ for $j=1,\dots, k$, where the
words $r_{f(t)}$ and $r_u$ were defined by formulas (\ref{rel10})
and (\ref{rel11}). Denote by
$r_{m+k+i}:=x_{m+1}^nx_ix_{m+1}^{-n}x_i^{-1}$ if
$$u_{m+k+i}= (0, \dots, 0,t^n-1,0,\dots,0)\in M_1$$  for
$i=1,\dots, m$ and denote by $\overline G$ the group given by
presentation
$$\overline G=\langle x_1,\dots,x_{m+1} \mid r_1,\dots,r_{2m+k}\rangle .$$

It follows from Lemma \ref{ftor} that the matrix
$\widetilde{\mathcal A}(G)$ (resp., $\widetilde{\mathcal
A}(\overline G)$), formed by the first $m$ columns of the
Alexander matrix $\mathcal A(G)$ (resp., $\mathcal A(\overline
G)$), coincides with the matrix $\mathcal U$ (resp.,
$\overline{\mathcal U}$) formed by the rows $u_1,\dots,u_{m+k}$
(resp., by $u_1,\dots,u_{2m+k}$). Therefore, by Proposition
\ref{kum}, the Alexander module $A_0(G)$ (resp., $A_0(\overline
G)$) coincides with $M=\Lambda^m/M_1$, where $M_1$ is generated by
the rows $u_1,\dots,u_{m+k}$ (resp., by $u_1,\dots,u_{2m+k}$).

Notice that $G$ (resp., $\overline G$) is an irreducible
$C$-group, since all $C$-generators $x_1,\dots,x_m$ are conjugated
to $x_{m+1}$. Moreover, $\overline G$ is a Hurwitz $C$-group.
Indeed, it follows from relations $r_{m+k+j}$, $j=1,\dots,m$, that
$x_{m+1}^n$ belongs to the center of $\overline G$. Since all
$x_i$ are conjugated to $x_{m+1}$, we have $x_i^n=x_{m+1}^n$ for
all $i=1,\dots,m$. Therefore the product $x_1^n\dots x_{m+1}^n$
also belongs to the center of $\overline G$ and $\overline G$
possesses a Hurwitz presentation
$$\begin{array}{ll}
\overline G=\langle x_1,\dots ,x_{n(m+1)}\mid & r_1,\dots
,r_{2m+k},
\\ & x_ix_{i+m+1}^{-1},\,\, i=1,\dots ,(n-1)(m+1), \\
& [x_i,(x_1\dots x_{n(m+1})],\, \, i=1,\dots ,n(m+1)\rangle .
\end{array}
$$

The following two lemmas complete the proof of Theorems
\ref{main1} and \ref{main2}.
\begin{lem} \label{kuz} {\rm (\cite{Kuz2})} The Alexander module
$A_0(G)=G^{\prime}/G^{\prime\prime}$ of an irreducible $C$-group
$G$ is a Noetherian $(t-1)$-invertible $\Lambda$-module.
\end{lem}
\proof For an irreducible $C$ group $G$ its commutator subgroup
$G^{\prime}$ coincides with the kernel of the $C$-epimorphism $\nu
:G\to \mathbb F_1$. By the Reidemeister -- Schreier method, if
$C$-generators $x_1,\dots,x_m$ generate $G$, then the elements
$a_{i,n}=x_m^nx_ix_m^{-(n+1)}$, $i=1,\dots,m-1$, $n\in \mathbb Z$,
generate $G^{\prime}$. Therefore
$A_0(G)=G^{\prime}/G^{\prime\prime}$ is generated by the images
$\overline a_{i,n}$ of the elements $a_{i,n}$ under the natural
epimorphism $G^{\prime}\to G^{\prime}/G^{\prime\prime}$. The
action of $t$ on $A_0(G)$ is defined by conjugation $a\mapsto
x_max_m^{-1}$ for $a\in G^{\prime}$. Therefore $t\overline
a_{i,n}=\overline a_{i,n+1}$. Thus $A_0(G)$ is generated over
$\Lambda$ by $\overline a_{1,0}\dots, \overline a_{m-1,0}$ and
hence it is a Noetherian $\Lambda$-module.

To show that $A_0(G)$ is a $(t-1)$-invertible $\Lambda$-module,
notice, first, that any element $g\in G$ can be written in the
form $g=x_m^ka$, where $a\in G^{\prime}$ and $k=\nu(g)$. Therefore
$G^{\prime}$ is generated by the elements of the form
$[x_m^na,x_m^kb]$, where $a,b\in G^{\prime}$, and hence $A_0(G)$
is generated by their images $\overline{[x_m^na,x_m^kb]}$. It is
easily to check that
\begin{equation} \label{u}
\begin{array}{ll}
[x_m^na,x_m^kb]=   &   [ x_m^{n},a]
(ax_m^{n+k}[b,a^{-1}]x_m^{-(n+k)}a^{-1})[a,x_m^{n+k}]\cdot
\\ & \cdot (x_m^{n+k}[b,x_m^{-n}]x_m^{-(n+k)}).
\end{array}
\end{equation}
It follows from (\ref{u}) that
\begin{equation} \label{u2}
\begin{array}{ll}
\overline{[x_m^na,x_m^kb]}= & (t^n-1)\overline a+
(1-t^{n+k})\overline a+t^{n+k}(1-t^{-n})\overline b= \\
& t^n(1-t^k)\overline a +t^k(t^n-1)\overline b= \\
& \displaystyle (t-1)(\sum_{i=0}^{n-1}t^{i+k}\overline b
-\sum_{i=0}^{k-1}t^{i+n}\overline a),
\end{array}
\end{equation}
since $ax_m^{n+k}[b,a^{-1}]x_m^{-(n+k)}a^{-1}\in
G^{\prime\prime}$. Now, it is easy to see that the multiplication
by $t-1$ is an epimorphism of $A_0(G)$, since the elements of the
form $\overline{[x_m^na,x_m^kb]}$ generate $A_0(G)$ over $\mathbb
Z$. To complete the proof, we apply Lemma \ref{crit0}. \qed

\begin{lem} \label{id} {\rm (\cite{Ku?})}
The Alexander module of a Hurwitz $C$-group of degree $m$ is a
Noetherian $(t-1)$-invertible $t$-unipotent $\Lambda$-module of
unipotence index $d$, where $d$ is a divisor of $m$.
\end{lem}
\proof If $G$ is a Hurwitz group of degree $m$, then it is
generated by $C$-generators $x_1,\dots , x_m$ such that the
product $x_1\dots x_m$ belongs to the center of $G$. By Lemma
\ref{kuz}, the Alexander module
$A_0(G)=G^{\prime}/G^{\prime\prime}$ in a Noetherian
$(t-1)$-invertible $\Lambda$-module. The multiplication by $t$ on
$A_0(G)$ is induced by conjugation $a\mapsto x_max_m^{-1}$ for
$a\in G'$. Since $\nu(x_m^m)=\nu (x_1\dots x_m)$, there is an
element $a_0\in G'$ such that $x_m^m=a_0\cdot x_1\dots x_m$ and
hence the conjugation by $x_m^m$ is an inner automorphism of $G'$.
Therefore the induced automorphism
$t^m$ of $G^{\prime}/G^{\prime\prime}$ is the identity. \qed \\

\subsection{Alexander modules of $C$-products of $C$-groups}
Let $G_1$ and $G_2$ be two irreducible $C$-groups and let $x\in
G_1$ (resp., $y\in G_2$) be one of $C$-generators of $G_1$ (resp.,
of $G_2$). Consider the amalgamated product $G_1*_{\{ x=y\}}G_2$.
If
\begin{equation} \label{equa1} \begin{array}{l} G_1=\langle x_1,\dots, x_n \mid \mathcal
R_1\rangle ,\\
G_2=\langle y_1,\dots, y_m \mid \mathcal R_2\rangle \end{array}
\end{equation}
are $C$-presentations of $G_1$ and $G_2$, where $x=x_n$ and
$y=y_m$, then $G_1*_{\{ x=y\}}G_2$ is given by $C$-presentation
\begin{equation} \label{equa2} \langle x_1,\dots, x_{n-1}, y_1,\dots,y_{m-1},z \mid
\widetilde{\mathcal R}_1\cup \overline{\mathcal R}_2 \rangle
\end{equation}
in which each relation  $\widetilde r_i\in\widetilde{\mathcal
R}_1$ (resp., $\overline r_i\in\overline{\mathcal R}_2$) is
obtained from the relation $r_i\in\mathcal R_1$ (resp., from
$r_i\in\mathcal R_2$) by substitution of $z$ instead of $x_n$
(resp., instead of $y_m$).

If $x'\in G_1$ and $y'\in G_2$ are two another $C$-generators of
these groups, then there are inner $C$-isomorphisms $f_i:G_i\to
G_i$ such that $f_1(x')=x$ and $f_2(y')=y$, since all
$C$-generators of an irreducible $C$-group are conjugated to each
other. Therefore there is an isomorphism $$f_1*f_2:G_1*_{\{
x'=y'\}}G_2\to G_1*_{\{ x=y\}}G_2,$$ that is, the group $G_1*_{\{
x=y\}}G_2$, up to a $C$-isomorphism, does not depend on the choice
of $C$-generators $x$ and $y$, so we denote it by $G_1*_CG_2$ and
call the {\it $C$-product} of irreducible $C$-groups $G_1$ and
$G_2$.
\begin{prop}
\label{cprod} If a $C$-group $G=G_1*_CG_2$ is the $C$-product of
irreducible $C$-groups $G_1$ and $G_2$, then its Alexander module
$A_0(G)$ is isomorphic to the direct sum of the Alexander modules
of $G_1$ and $G_2$,
$$A_0(G)=A_0(G_1)\oplus A_0(G_2).$$
\end{prop}
\proof This proposition is a simple consequence of Proposition
\ref{kum}. Indeed, if $G_1$ and $G_2$ are given by presentation
(\ref{equa1}), then, by Proposition \ref{kum}, the Alexander
module $A_0(G)$ of the $C$-group $G=G_1*_CG_2$, given by
presentation (\ref{equa2}), is isomorphic to a factor module
$\Lambda^{n+m-1}/M(G)$, where the submodule $M(G)$ of
$\Lambda^{n+m-1}$ is generated by the rows of the matrix
$$
\overline{\mathcal A}=\left( \begin{array}{cc}
\overline{\mathcal A}_1 & 0 \\
0 &  \overline{\mathcal A}_2 \end{array} \right) , $$ where
$\overline{\mathcal A}_1$ (resp., $\overline{\mathcal A}_2$) is
the matrix formed by the first $n-1$ (resp., $m-1$) columns of the
matrix $\mathcal A(G_1)$ (resp., $\mathcal A(G_2)$). Now, it is
easy to see that $A_0(G)=A_0(G_1)\oplus A_0(G_2)$. \qed \\

Let \begin{equation} \label{pres-3} G=<x_1,\dots,x_m\, \mid\,
r_1,\dots,r_n> \end{equation} be a $C$-presentation of a $C$-group
$G$. The number $d_P=m-n$ is called the {\it $C$-deficiency} of
presentation (\ref{pres-3}) and $d_G=\min d_P$, where the minimum
is taken over all $C$-presentation of a $C$-group $G$, is called
the {\it $C$-deficiency} of the group $G$. Obviously, for a
$C$-group consisting of $k$ connected component, its
$C$-deficiency $d_G\leq k$ and, in particular, if $G$ is an
irreducible $C$-group, then $d_G\leq 1$.
\begin{lem} \label{defici-prod} Let $G=G_1*_CG_2$ be the
$C$-product of irreducible $C$-groups $G_1$ and $G_2$. Then
$$d_G\geq d_{G_1}+d_{G_2}-1.$$

In particular, if $d_{G_1}=d_{G_2}=1$, then $d_G=1$.
\end{lem}
\proof It follows from formula (\ref{equa2}).

\subsection{Presentation graphs of $C$-groups}
Let us associate a {\it presentation graph} $\Gamma_P$  to each
$C$-presentation (\ref{pres-3}) as follows. The vertices  of the
graph $\Gamma_P$ are labeled by the generators from presentation
(\ref{pres-3}) (and, in particular, they are in one to one
correspondence with these generators); the edges of $\Gamma_P$ are
in one to to one correspondence with the relations $r_j$ of the
presentation (\ref{pres-3}) and if
$r_j:=w_j^{-1}(x_1,\dots,x_m)x_{i_1}w_j(x_{1},\dots,x_m)x_{i_2}^{-1}$,
then the corresponding edge connects the vertices  $x_{i_1}$ and
$x_{i_2}$.

Obviously, the $C$-deficiency
$$d_P=\dim H_0(\Gamma_P,\mathbb
R)-\dim H_1(\Gamma_P,\mathbb R).$$
Therefore {\it for an
irreducible $C$-group $G$ its $C$-deficiency $d_G=1$ if and only
if $G$ possesses a $C$-presentation whose presentation graph
$\Gamma_P$ is a tree.}

A $C$-presentation
\begin{equation} \label{pres-2} G=<x_1,\dots,x_m\, \mid\,
r_1,\dots,r_n>\end{equation} is said to be {\it simple} if each
relation $r_j$ in (\ref{pres-2}) is of the form:
$$r_j:=x_{i_3}^{-1}x_{i_1}x_{i_3}x_{i_2}^{-1},$$
for some $i_1,i_2,i_3\in \{ 1,\dots, m\}$ (that is,
$x_{i_2}=x_{i_3}^{-1}x_{i_1}x_{i_3}$).

\begin{rem} \label{rem?} If presentations {\rm (\ref{equa1}\rm )} of irreducible $C$-groups
$G_1$ and $G_2$ are simple, then presentation {\rm (\ref{equa2}\rm
)} of $G=G_1*_CG_2$ is also simple and the presentation graph
$\Gamma_P$ of presentation {\rm (\ref{equa2}\rm )} is the bouquet
$\Gamma_P=\Gamma_{P_1}\bigvee_{z=x_n=y_m}\Gamma_{P_2}$ of the
presentation graphs $\Gamma_{P_1}$ and $\Gamma_{P_2}$ of
presentations {\rm (\ref{equa1}\rm )}. In particular, if
$\Gamma_{P_1}$ and $\Gamma_{P_2}$ are trees, then the presentation
graph $\Gamma_P$ is also a tree.
\end{rem}
\begin{lem} \label{simpres} Any $C$-group possesses a simple
$C$-presentation with $C$-de\-ficien\-cy $d_P=d_G$. \end{lem}

\proof Let $G$ be given by $C$-presentation of $C$-deficiency
$d_P=d_G$ and $r:=w^{-1}x_iwx_j^{-1}$ is one of its relations
(that is, $w^{-1}x_iw=x_j$), where $w=x_{i_1}^{\varepsilon_1}\dots
x_{i_k}^{\varepsilon_k}$ is a word in $\mathbb F_m$ and
$\varepsilon_l=\pm 1$, then we can add $k-1$ new generators
$x_{m+1}$,$\dots,$ $x_{m+k-1}$ and replace  the relation $r$ by
$k$ relations:
$$\begin{array}{lcl}
x_{m+1}& = & x_{i_1}^{-\varepsilon_1}x_ix_{i_1}^{\varepsilon_1}, \\
x_{m+2} & = &
x_{i_2}^{-\varepsilon_2}x_{m+1}x_{i_2}^{\varepsilon_2},
\\\dots \dots & \dots & \dots \dots \dots \dots \\
x_{m+k-1} & = &
x_{i_{k-1}}^{-\varepsilon_{k-1}}x_{m+k-2}x_{i_{k-1}}^{\varepsilon_{k-1}},
\\ x_{j} & = & x_{i_k}^{-\varepsilon_k}x_{m+k-1}x_{i_k}^{\varepsilon_k}.
\end{array}
$$
Obviously, we obtain a new $C$-presentation of the same
$C$-deficiency which defines the same $C$-group $G$. \qed

\subsection{The Alexander modules of $C$-groups possessing
$C$-presentations whose presentation graphs are trees} By Lemma
\ref{simpres}, an irreducible $C$-group $G$ possesses a simple
$C$-presentation whose presentation graph is a tree if and only if
its $C$-deficiency $d_G=1$.

\begin{prop}
\label{treepres} If $M=\bigoplus_{i=1}^mM_i$ is the direct sum of
bi-principle $(t-1)$-invertible $\Lambda$-modules
$M_i=\Lambda/\langle f_i(t)\rangle$, then there is an irreducible
$C$-group $G$ such that $A_0(G)\simeq M$ and such that its
$C$-deficiency $d_G=1$.
\end{prop}
\proof Note that the $C$-deficiency of a $C$-group given by
presentation
\begin{equation} \label{sineq} G=\langle x_1,x_2\mid
wx_1 w^{-1}x_{2}^{-1}\rangle,
\end{equation}
where $w=w(x_1,x_2)$ is a word in letters $x_1$, $x_2$ and their
inverses, is equal to $1$. Applying Proposition \ref{kum}, we see
that the Alexander module $A_0(G)$ of a $C$-group $G$, given by
presentation (\ref{sineq}), is a bi-principle $(t-1)$-invertible
$\Lambda$-module.

Conversely, it was shown in the proof of Theorem \ref{main11} that
any bi-principle $(t-1)$-invertible $\Lambda$-module
$M=\Lambda/\langle f(t)\rangle$ is the Alexander module of some
irreducible $C$-group given by presentation (\ref{sineq}). To
complete the proof we apply Proposition \ref{cprod} and Remark
\ref{rem?}.\qed

\begin{cor} \label{sfera}
Let $M=\bigoplus_{i=1}^mM_i$ is a direct sum of bi-principle
$(t-1)$-invertible $\Lambda$-modules $M_i=\Lambda/\langle
f_i(t)\rangle$. Then for each $n\geq 2$ there is a knotted sphere
$S^n\subset S^{n+2}$ such that the Alexander module
$$A_0(\pi_1(S^{n+2}\setminus S^n))\simeq M.$$

In particular, a polynomial $f(t)\in \mathbb Z[t]$ is the
Alexander polynomial $\Delta (t)$ of some knotted sphere
$S^n\subset S^{n+2}$ if and only if $f(1)=\pm 1$ and, moreover,
the Jordan blocks of the Jordan canonical form of the matrix of
the automorphism $h_{\mathbb C}$ acting on $A_0(S^n)\otimes
\mathbb C$ can be of arbitrary size.
\end{cor}
\proof In \cite{Ku1}, it was shown that if an irreducible
$C$-group is given by a simple $C$-presentation which presentation
graph is a tree, then for each $n\geq 2$ there is a knotted sphere
$S^n\subset S^{n+2}$ such that $\pi_1(S^{n+2}\setminus S^n)\simeq
G$. \qed
\begin{prop}
Let $G$ be an irreducible $C$-group of $C$-deficiency $d_G=1$.
Then its Alexander module $A_0(G)$ has not non-zero $\mathbb
Z$-torsion elements.
\end{prop}
\proof Let
\begin{equation} \label{pres-5} G=<x_1,\dots,x_m\, \mid\,
r_1,\dots,r_{m-1}>\end{equation} be a $C$-presentation of $G$. By
Proposition \ref{kum}, its Alexander module $A_0(G)$ is isomorphic
to a factor module $\Lambda^{m-1}/M(G)$, where the submodule
$M(G)$ of $\Lambda^{m-1}$ is generated by the rows of the matrix
$\overline{\mathcal A}$ formed by the first $m-1$ columns of the
Alexander matrix $\mathcal A(G)$ of the group $G$ given by
presentation (\ref{pres-5}). The size of the matrix
$\overline{\mathcal A}$ is $(m-1)\times(m-1)$.
\begin{lem} \label{determ} The determinant
$\Delta(t)=\det \overline{\mathcal A}$ satisfies the following
property: $\Delta(1)=\pm 1$.
\end{lem}
\proof It coincides with the similar statement for knot groups
(see the proof, for example, in \cite{C-F}). \qed

Denote by $\mathcal A_j$ the rows of the matrix
$\overline{\mathcal A}$, $j=1,\dots, m-1$. The module $A_0(G)$ has
a non-zero $\mathbb Z$-torsion element if and only if there is a
vector $u=(f_1(t),\dots,f_{m-1}(t))$ such that $u\not\in M(G)$ and
$ku\in M(G)$ for some $k\in\mathbb N$. Assume that there is a such
vector $u$. Then there are $g_j(t)\in \Lambda$ such that $ku=\sum
g_j(t)\mathcal A_j$, where for some $g_j(t)$ one of its
coefficients is not divisible by $k$.

Without loss of generality, we can assume that all $f_i(t)$ and
$g_j(t)$ belong to $\mathbb Z[t]$. By Cramer's theorem,
$$g_j(t)=\frac{\Delta_j(t)}{\Delta(t)},$$
where $\Delta_j(t)$ is the determinant of  the matrix obtained
from $\overline{\mathcal A}$ by substitution $ku$ instead of the
row $\mathcal A_j$. Therefore the coefficients of all polynomials
$\frac{\Delta_j(t)}{\Delta(t)}$ are divisible by $k$. A
contradiction. \qed
\begin{rem} If $G$ is an irreducible $C$-group given by presentation of $C$-deficiency
$d_P=d_G=1$, then the determinant $\Delta(t)=\det
\overline{\mathcal A}$ of the matrix $\overline{\mathcal A}$,
obtained from the Alexander matrix $\mathcal A$ after deleting its
last column, coincides with the Alexander polynomial $\Delta_G(t)$
of the group $G$.
\end{rem}

\subsection{Finitely $\mathbb Z$-generated Alexander modules of
irreducible $C$-groups}
\begin{thm} \label{finAmod} Let $G$ be an irreducible $C$-group.
The Alexander module $A_0(G)$ is finitely generated over $\mathbb
Z$ if and only if the leading coefficient $a_n$ and the constant
coefficient $a_0$ of the Alexander polynomial
$\Delta_G(t)=\sum_{i=0}^na_it^i$ of  $G$ are equal to $\pm 1$.
\end{thm}
\proof By Theorem \ref{main11}, $A_0(G)$ is a Noetherian
$(t-1)$-invertibele $\Lambda$-module. Let $A_0(G)_{fin}$ be the
$Z$-torsion submodule of the Alexander module $A_0(G)$. By Theorem
\ref{fin=fin}, $A_{0}(G)_{fin}$ is finitely generated over
$\mathbb Z$.

Consider the quotient module $M=A_0(G)/A_{0}(G)_{fin}$. It is free
from $\mathbb Z$-torsion. Therefore there is a natural embedding
$M\hookrightarrow M_{\mathbb Q}=M\otimes \mathbb Q$. We have
$\dim_{\mathbb Q} M_{\mathbb Q}<\infty$, since $M$ is a Noetherian
$\Lambda$-torsion module.

Denote by $h_{\mathbb Q}$ an automorphism of $M_{\mathbb Q}$
induced by the multiplication by $t$. Then, by definition,
$\Delta_G(t)=a\det(h_{\mathbb Q}-t\text{Id})$, where $a\in \mathbb
N$ is the smallest number such that $a\det (h_{\mathbb
Q}-t\text{Id})\in \mathbb Z[t]$.

If the Alexander module $A_0(G)$ is finitely generated over
$\mathbb Z$, then $M$ is a free finitely generated $\mathbb
Z$-module. Denote by $h$ an automorphism of $M$ induced by
multiplication by $t$. We have $\det h=\pm 1$ and
$$\det(h-t\text{Id})=\det(h_{\mathbb Q}-t\text{Id})\in \mathbb Z[t].$$
Therefore $\Delta_G(t)=\det(h-t\text{Id})$ and its leading
coefficient $a_n=(-1)^n$, where $n=\text{rk}\, M$, and $a_0=\det
h=\pm 1$.

Let the leading coefficient $a_n$ and the constant coefficient
$a_0$ of the Alexander polynomial $\Delta_G(t)$ of $G$ be equal to
$\pm 1$. By Cayley-Hamilton's Theorem, $\Delta_G(t)\in
\text{Ann}(M_{\mathbb Q})$. Therefore $\Delta_G(t)\in
\text{Ann}(M)$ and $M$ is finitely generated over $\mathbb Z$ by
Proposition \ref{finitenes}. \qed
\begin{rem} Let an irreducible $C$-group $G$ is given by
$C$-presentation $G=\langle x_1,\dots,x_m\mid
r_1,\dots,r_n\rangle$ and $\mathcal A(G)$ its Alexander matrix.
Then the Alexander polynomial $\Delta_G(t)$ coincides (up to
multiplication by $\pm t^k$) with the greatest common divisor of
the determinants of all $(m-1)\times(m-1)$ submatrices $\mathcal
A_{m-1}$ of the matrix $\mathcal A(G)$.
\end{rem}

\subsection{Alexander modules of some irreducible
$C$-groups}
In the end of this section, we compute the Alexander
modules for some irreducible $C$-groups.
\begin{ex} \label{ex4} The Alexander module $A_0(\text{Br}_{m+1})$ of the braid group
$\text{Br}_{m+1}$ is trivial if $m\geq 4$ {\rm (}or $m=1${\rm )}
and isomorphic to $\Lambda/\langle t^2-t+1\rangle$  for $m=2$ and
$3$.
\end{ex}
This statement is well known, but for completeness, we give a
proof. \proof The braid group $\text{Br}_{m+1}$ is given by
presentation
$$\begin{array}{lcl} \text{Br}_{m+1}=
 \langle x_1,\dots, x_{m} \mid & [x_i,x_j] & \text{for}\, \, |i-j|\geq 2, \\
 & x_ix_{i+1}x_ix_{i+1}^{-1}x_i^{-1}x_{i+1}^{-1}\quad & \text{for}\, \, i=1,\dots,m-1 \rangle .
 \end{array}$$
Notice that it is a $C$-presentation of an irreducible $C$-group.

By Proposition \ref{kum}, to calculate $A_0(\text{Br}_{m+1})$ we
should calculate the matrix $\overline{\mathcal
A}(\text{Br}_{m+1})$.

The relations $[x_m,x_{i}]$, $i=1,\dots,m-2$, give the rows
\begin{equation} \label{rows} (0,\dots,0,(t-1),0\dots,0),
\end{equation} where $t-1$ stands on the $i$-th place for $i=1,\dots,m-2$,
and if $m\geq 4$, then the relation $[x_{m-1},x_{1}]$ gives the
row
\begin{equation} \label{rows1} (t-1,0,\dots,0,1-t).
\end{equation}
If $m\geq 4$, then the rows from (\ref{rows}) and row
(\ref{rows1}) generate submodule $(t-1)\Lambda^{m-1}$ of the
module $\Lambda^{m-1}$. On the other hand, these rows belong to
the module $M(\text{Br}_{m+1})$. Therefore
$A_0(\text{Br}_{m+1})=0$, since $A_0(\text{Br}_{m+1})\simeq
\Lambda^{m-1}/M(\text{Br}_{m+1})$ is a $(t-1)$-invertible
$\Lambda$-module and $(t-1)\Lambda^{m-1}\subset
M(\text{Br}_{m+1})$.

If $m=2$, then we have the only one relation in the presentation
of $\text{Br}_{3}$, namely,
$$r:=x_1x_{2}x_1x_{2}^{-1}x_1^{-1}x_{2}^{-1}$$
We have $\nu_*(\frac{\partial r}{\partial x_1})=1+t^2-t$ and
therefore $A_0(\text{Br}_{3})\simeq \Lambda/\langle
t^2-t+1\rangle$.

If $m=3$, then we have the only three relations in the
presentation of $\text{Br}_{4}$, namely,
$$\begin{array}{ll} r_1:= & x_1x_{2}x_1x_{2}^{-1}x_1^{-1}x_{2}^{-1}, \\
 r_2:= & x_2x_{3}x_2x_{3}^{-1}x_2^{-1}x_{3}^{-1}, \\
  r_3:= & x_1x_{3}x_{1}^{-1}x_3^{-1}.
\end{array}$$
We have $$\begin{array}{l} \nu_*(\frac{\partial r_1}{\partial
x_1})=-\nu_*(\frac{\partial r_1}{\partial
x_2})=\nu_*(\frac{\partial r_2}{\partial x_2})=t^2-t+1 \\
\nu_*(\frac{\partial r_3}{\partial x_1})=1-t,
\end{array}$$
Therefore $M(\text{Br}_{3})\subset \Lambda^2$ is generated by
vectors $$v_1=(t^2-t+1,-(t^2-t+1)),\quad v_2=(0,t^2-t+1), \quad
v_3=(1-t,0),$$ and hence $A_0(\text{Br}_{3})\simeq \Lambda/\langle
t^2-t+1\rangle$.\qed
\begin{ex}\label{ex5}
The Alexander module of a $C$-group
$$G_m=\langle x_1,x_2 \mid
(x_1^{-1}x_2)^mx_1(x_1^{-1}x_2)^{-m}x_2^{-1}\rangle ,$$ $m\in
\mathbb N$, is isomorphic to $A_0(G)\simeq \Lambda/\langle
(m+1)t-m\rangle$.
\end{ex}
These irreducible $C$-groups are interesting, since they are
non-Hopfian if $m\geq 2$ and therefore they are not residually
finite. (The group $G_m$ is isomorphic to Baumslag -- Solitar
group (see \cite{B-S}) $\langle a,x_1\mid
x_1^{-1}a^mx_1a^{-(m+1)}\rangle$ if we put $x_2=x_1a$.) Note also
that each of these groups can be realized as $\pi_1(S^4\setminus
S^2)$ for some knotted sphere $S^2\subset S^4$.

\proof Straightforward calculation gives
$$\nu_*(\frac{\partial r}{\partial
x_1})=-mt^{-1}+m+1, $$ where
$r:=(x_1^{-1}x_2)^mx_1(x_1^{-1}x_2)^{-m}x_2^{-1}$. Therefore the
Alexander module $$A_0(G)\simeq \Lambda/\langle (m+1)t-m\rangle.
\qquad \qquad \qed $$

\section{First homology groups of cyclic coverings}
\subsection{Proof of Theorems \ref{kmain1} and \ref{hmain2}}
Theorems \ref{kmain1} and \ref{hmain2} will be proved
simultaneously.

In the notations from Introduction, we denote by $X$ either the
sphere $S^{n+2}$ (\it Case I}) or $\mathbb C\mathbb P^2$ ({\it
Case II}), and by $X'$ respectively either the complement of a
knotted $n$-manifold $V$ in $S^{n+2}$ or the complement of the
union of an irreducible Hurwitz curve $H$ and a line\, "at
infinity" $L$ in $\mathbb C\mathbb P^2$. Recall that the
fundamental group $G=\pi_1(X')$ is an irreducible $C$-group.

Consider the infinite cyclic covering $f =f_{\infty }:X_{\infty }
\to X'$ corresponding to the $C$-epimorphism $\nu :G\to \mathbb
F_1$ with $\ker \nu=G'$. Let $h\in \text{Deck}(X_{\infty }
/X')\simeq {\mathbb F}_{1}$ be a covering transformation
corresponding to the $C$-generator $x\in {\mathbb F}_{1}$. We say
that $h$ is the {\it monodromy} respectively of the knotted
manifold $V$ and of the Hurwitz curve $H$. The space $X'$ will be
considered as the quotient space $X'=X_{\infty }/{\mathbb F}_{1}$.
In such a situation Milnor \cite{Mi} considered an exact sequence
of chain complexes
$$ 0\to C_{\cdot }(X_{\infty })\buildrel{h-id}\over\longrightarrow
 C_{\cdot }(X_{\infty })\buildrel{f_{*}}\over\longrightarrow
 C_{\cdot }(X')\to  0 $$
which gives an exact sequence of homology groups with integer
coefficients:
\begin{equation} \label{chain} \ldots \to H_{1}(X_{\infty
})\stackrel{t-id}\longrightarrow H_{1}(X_{\infty })\stackrel{f
_*}\longrightarrow H_{1}(X')\stackrel{\partial}\longrightarrow
H_{0}(X_{\infty })\rightarrow  0 ,\end{equation} where $t=h_{*}$.

The action $h_*$ (resp., $h_{k*}$) defines on $H_{1}(X_{\infty
})\simeq G^{\prime}/G^{\prime\prime}$ a structure of
$\Lambda$-module such that sequence (\ref{chain}) is an exact
sequence of $\Lambda$-modules (so that $H_{1}(X_{\infty })$ is the
Alexander module of the $C$-group $G$). The action of
$t\in\Lambda$ on $H_{0}(X_{\infty })\simeq \mathbb Z$ is trivial,
that is, $t$ is the identity automorphism of $H_{0}(X_{\infty })$.

If $\langle h^k\rangle\subset {\mathbb F}_{1}$ is an infinite
cyclic group generated by $h^{k}$, then $X'_k= X_{\infty }/\langle
h^k\rangle$ and $X'= X'_k/\mu _{k}$, where $\mu _{k}={\mathbb
F}_{1}/\langle h^k\rangle$ is the cyclic group of order $k$.
Denote by $h_{k}$ an automorphism of $X'_k$ induced by the
monodromy $h$. Then $h_{k}$ is a generator of the covering
transformation group $\text{Deck}(X'_{k} /X')=\mu _{k}$ acting on
$X'_k$.

It is easy to see that in Case I the manifold $X'_k$ can be
embedded to the compact smooth manifold $X_k$ satisfying the
following properties:
\begin{itemize}
\item[($i$)] the action of $h_k$ on $X'_k$ and the map
$f'_k:X'_k\to X'$ are continued to an action (denote it again by
$h_k$) on $X_k$ and to a smooth map
$$f_k:X_k\to X\simeq X_k/\{ h_k\},$$
\item[($ii$)] the set of fixed points of $h_k$ coincides with
$f_k^{-1}(V)=\overline V$ and the restriction $f_{k\mid \overline
V}:\overline V\to V$ of $f_k$ to $\overline V$ is a smooth
isomorphism.
\end{itemize}

In Case II (in the notations of the proof of Theorem 4.1 in
\cite{Gr-Ku}), the covering $f'_k$ can be extended to a map
${\widetilde f}_{k\, {\rm norm}}:{\widetilde X}_{k\, {\rm norm}}
\to X$ branched along $H$ and, maybe, along $L$.  Let $\sigma
:\overline X_k\to \widetilde X_{k,\, \text{norm}}$ be a resolution
of the singularities, $E=\sigma^{-1}(\text{Sing}\, \widetilde
X_{k,\, \text{norm}})$, and $\overline f_k=\widetilde f_{k,\,
\text{norm}}\circ \sigma$. Denote by $R=\widetilde f_{k,\,
\text{norm}}^{-1}(H)$ and $R_{\infty}=\widetilde f_{k,\,
\text{norm}}^{-1}(L)$. Note that the restriction of $\widetilde
f_{k,\, \text{norm}}$ to $R$ is one-to-one and the restriction of
$\widetilde f_{k,\, \text{norm}}$ to $R_{\infty}$ is a
$k_0$-sheeted cyclic covering, where $k_0=\text{GCD}(k,d)$ and the
ramification index of $\widetilde f_{k,\, \text{norm}}$ along
$R_{\infty}$ is equal to $k_{\infty}=\frac{k}{k_0}$. As in the
algebraic case, it is easy to show that $R_{\infty}$ is
irreducible. Denote by $\overline R=\sigma^{-1}(R)$ the proper
transform of $R$. Note that $k_0$ is a divisor of $m$. Put
$m_0=\frac{m}{k_0}$, we have $m_0\in \mathbb N$.

Denote by $X_k=\overline X_k\setminus E$. We have two embeddings
$i_{k}:X'_k\hookrightarrow X_k$ and $j_{k}:X_k\hookrightarrow
\overline X_k$.

In both cases , the action of $h_k$ on $X_k$ induces on
$H_1(X_k,\mathbb Z)$ (resp., on $H_{1}(X'_k,\mathbb Z)$) a
structure of $\Lambda$-module such that the homomorphism
$$i_{k*}:H_1(X'_k,\mathbb Z) \to H_1(X_k,\mathbb Z),$$ induced by
the embedding $i:X'_k\hookrightarrow X_k$, is a
$\Lambda$-homomorphism. Obviously, the homomorphism $i_{k*}$ is an
epimorphism.

In Case I, let $S\subset X_k$ be a germ of a smooth surface
meeting transversally $\overline V$ at $p\in \overline V$ and let
$\bar \gamma\subset S$ be a circle of small radius with center at
$p$. Then $\ker i_{k*}$ is generated by the homology class
$[\bar\gamma ]\in H_1(X'_k,\mathbb Z)$ containing the cycle
$\bar\gamma$, since $\overline V$ is a smooth connected
codimension two submanifold of $X_k$.

It is obvious, that $t([\bar\gamma])=[\bar\gamma]$, where
$t=h_{k*}$, and
$$f_{k*}([\bar\gamma])=\pm k[\gamma]\in H_1(X',\mathbb Z)\simeq
\mathbb Z,$$ where $[\gamma]$ is a generator of $H_1(X',\mathbb
Z)$ represented by  a simple loop $\gamma$ around $V$.

In Case II, let $S\subset X_k$ be a germ of a smooth surface
meeting transversally $R$ at $p\in R$ and let $\bar \gamma\subset
S$ be a circle of small radius with center at $p$. Evidently, the
homology class $[\bar \gamma]\in H_1(X'_k,\mathbb Z)$ is invariant
under the multiplication by $t$ and
$f_{k*}([\bar\gamma])=k[\gamma]$, where $[\gamma]$ is a generator
of $H_1(\mathbb C\mathbb P^2\setminus(H\cup L),\mathbb Z)\simeq
\mathbb Z$.

Similarly, let a complex line $L_1\subset \mathbb C\mathbb P^2$
meet $L$ transversely  at $q\in L\setminus H$ and
$\gamma_{\infty}$ be a simple small loop around $L$ lying in
$L_1$. Then $f_k^{-1}(\gamma_{\infty})$ splits into the disjoint
union of $k_0$ simple loops $\bar \gamma_{\infty,i}$,
$i=1,\dots,k_0$. Since $R_{\infty}$ is irreducible, each two loops
$\bar \gamma_{\infty,i}$ and $\bar \gamma_{\infty,j}$ belong to
the same homology class of $H_1(X'_k,\mathbb Z)$ (denote it by
$[\bar \gamma_{\infty}]$). It is easy to see that $t(\bar
\gamma_{\infty,i})=\bar \gamma_{\infty,i+1}$. Therefore the
homology class $[\bar\gamma_{\infty}]\in H_1(X'_n\mathbb Z)$ is
invariant under the multiplication by $t$. Note also that
$f_{k*}([\bar\gamma_{\infty}])=k_{\infty}m[\gamma]=km_0[\gamma]$,
since $[\gamma_{\infty}]=m[\gamma]$.


\begin{lem} \label{splits}
The $\Lambda$-module $H_{1}(X'_k,\mathbb Z)$ is isomorphic to
$$A_k(G)\oplus H_1(X'_k)_1\simeq A_k(G)\oplus \mathbb Z,$$
where $A_k(G)$ is the $k$-th derived Alexander module of $C$-group
$G$ and $$H_1(X'_k)_1=\{ h\in H_{1}(X'_k,\mathbb Z) \mid
(t-1)h=0\}.$$
\end{lem}
\proof We apply the sequence
\begin{equation} \label{hexact}
\ldots \to H_{1}(X_{\infty },\mathbb
Z)\stackrel{t^k-id}\longrightarrow H_{1}(X_{\infty },\mathbb
Z)\stackrel{g_{k,*}}\longrightarrow H_{1}(X'_k,\mathbb
Z)\stackrel{\partial}\longrightarrow H_{0}(X_{\infty },\mathbb
Z)\rightarrow 0 \end{equation}
 constructed in the same way as (\ref{chain})
to the infinite cyclic covering $g_k=g_{\infty ,k}:X_{\infty }\to
X'_k$, to analyze the group $H_{1}(X'_{k},\mathbb Z)$.

By (\ref{hexact}), we have the short exact sequence
\begin{equation} \label{hexact1}
0\to H_{1}(X_{\infty })/(t^k-1)H_{1}(X_{\infty
})\stackrel{g_{k,*}}\longrightarrow
H_{1}(X'_k)\stackrel{\partial}\longrightarrow H_{0}(X_{\infty
})\rightarrow  0 \end{equation} which is a sequence of
$\Lambda$-homomorphisms.

Denote by $M_1=\ker \partial =\text{im} g_{k,*}\simeq
H_{1}(X_{\infty })/(t^k-1)H_{1}(X_{\infty })$ and by $M_2=
H_{1}(X'_k)_1$.

We have $H_{0}(X_{\infty },\mathbb Z)\simeq \mathbb Z$. Let us
choose a generator $u\in H_{0}(X_{\infty },\mathbb Z)$ and let
$v_1\in H_{1}(X'_k,\mathbb Z)$ be an element such that $\partial
(v_1)=u$. Then $(t-1)v_1\in\ker \partial$, since $H_{0}(X_{\infty
},\mathbb Z)$ is a trivial $\Lambda$-module and $\partial$ is a
$\Lambda$-homomorphism. We fix a such $v_1$.

By Theorems \ref{main1} and \ref{main2}, $H_{1}(X_{\infty
},\mathbb Z)=A_0(G)$ is a Noetherian $(t-1)$-invertible
$\Lambda$-module. Therefore, by Proposition \ref{cor-krit3},
$$M_1\simeq H_{1}(X_{\infty })/(t^k-1)H_{1}(X_{\infty
})=A_k(G)$$ is also a Noetherian $(t-1)$-invertible
$\Lambda$-module and, by Theorem \ref{cor-krit4}, there is a
polynomial $g_1(t)\in \text{Ann}(M_1)$ such that $g_1(1)=1$. We
fix a such polynomial $g_1(t)$.

Consider the element $\overline v_1=g_1(t)v_1$. We have $\partial
(\overline v_1)=g_1(1)u=u$ and hence
$$(t-1)\overline v_1=(t-1)g_1(t)v_1=g_1(t)(t-1)v_1=0,$$
since $(t-1)v_1\in M_1$. Therefore $\overline v_1\in M_2$.

Note that $M_1\cap M_2=0$, since $M_1$ is $(t-1)$-invertible.
Therefore  $\partial$ maps $M_2$ isomorphically onto
$H_{0}(X_{\infty },\mathbb Z)$, that is, exact sequence
(\ref{hexact1}) splits and hence $H_{1}(X'_k,\mathbb Z)\simeq
M_1\oplus M_2$. \qed

\begin{lem} \label{image} For
$f_{k*}:H_{1}(X'_k,\mathbb Z)\longrightarrow H_{0}(X',\mathbb Z)$
we have
\begin{itemize}
\item[($i$)] $\ker f_{k,*}=A_k(G)\subset H_{1}(X'_k,\mathbb Z)$,
\item[($ii$)] $\text{im}\, f_{k,*}= k\mathbb Z\subset \mathbb
Z\simeq H_{1}(X',\mathbb Z)$ and the restriction of $f_{k*}$ to
$H_1(X'_k)_1$ is an isomorphism of $H_1(X'_k)_1$ with its image.
\end{itemize}
\end{lem}
\proof The group $H_{1}(X',\mathbb Z)$ is isomorphic to
$G/G'\simeq \mathbb Z$. Similarly, the group $H_{1}(X'_k,\mathbb
Z)$ is isomorphic to $G_k/G'_k$, where $G_k=\ker \nu_k$,
$$\nu_k=\mod_k\circ\nu: G\to \mu_k=\mathbb
Z/\langle h^k\rangle,$$ and $f_{k*}:H_{1}(X'_k,\mathbb
Z)\rightarrow H_{1}(X',\mathbb Z)$ coincides with the homomorphism
$$i_{k*}:G_k/G'_k\to G/G'$$ induced by the embedding
$i_k:G_k\hookrightarrow G$.

Let the $C$-group $G$ be given by $C$-presentation (\ref{pres}).
To describe $\ker i_{k*}$ and $\text{im}\, i_{k*}$, let us
consider again the two-dimensional complex $K$ described in
section \ref{sec3}. The complex $K$ has a single vertex $x_0$, its
one dimensional skeleton is a bouquet of oriented circles $s_j$,
$1\leq j\leq m$, corresponding to the $C$-generators of $G$ from
presentation (\ref{pres}), and $K\setminus (\cup
s_i)=\bigsqcup_{j=1}^l\stackrel{\circ}D_j$ is a disjoint union of
open discs, where each disc $D_j$ corresponds to the relation
$r_i$ from presentation (\ref{pres}) (we denote here by $l$ the
number of relations $r_i$ in presentation  (\ref{pres})).

The embedding $i_k:G_k\hookrightarrow G$ defines an un-ramified
covering $f_k: K_k\to K$, where $K_k$ is a two-dimensional complex
consisting of $k$ vertices $p_1,\dots, p_k$, $f_k(p_j)=x_0$; the
preimage $f^{-1}(s_j)=\bigsqcup_{s=1}^k \overline s_{j,s}$ is the
disjoint union of $k$ edges $\overline s_{j,s}$, $1\leq s\leq k$;
and the preimage $f^{-1}(\stackrel{\circ}D_j)=\bigsqcup_{s=1}^k
\stackrel{\circ}{\overline D}_{j,s}$ is also the disjoint union of
$k$ open discs $\stackrel{\circ}{\overline D}_{j,s}$, $1\leq s\leq
k$.

Let $h_{k}$ be a generator of the covering transformation group
$\text{Deck}(K_{k} /K)=\mu _{k}$ acting on $K_k$. The
homeomorphism $h_k$ induces an action $h_{k*}$ on the chain
complex $C_{\cdot }(K_k)$ and an action $t$ on $H_i(K_k,\mathbb
Z)$ so that this action defines on $H_i(K_k,\mathbb Z)$ a
structure of $\Lambda$-module. It is easy to see that this
structure on $H_1(K_k,\mathbb Z)$ coincides with one on
$H_1(X'_k,\mathbb Z)$ defined above if we identify
$H_1(K_k,\mathbb Z)$ and $H_1(X'_k,\mathbb Z)$ by means of
isomorphisms $H_1(K_k,\mathbb Z)\simeq G_k/G'_k$ and
$H_1(X'_k,\mathbb Z)\simeq G_k/G'_k$.

Consider the sequence of chain complexes
$$C_{\cdot }(K_k)\buildrel{h_{k*}-id}\over\longrightarrow
 C_{\cdot }(K_k)\buildrel{f_{k*}}\over\longrightarrow
 C_{\cdot }(K)\to  0 .$$
It is easy to see that $\text{im}\, (h_{k*}-id)=\ker f_{k*}$ and
$$\ker (h_{k*}-id) =(\sum_{j=0}^{k-1}h_{k*}^j)C_{\cdot}(K_k).$$
Now the proof of Lemma  \ref{image} follows from the exact
sequence
\begin{equation} \label{chain2}
\begin{array}{l}
\ldots \to H_{1}(C_{\cdot}(K_k/\ker
(h_{k*}-id))\stackrel{t^k-1}\longrightarrow H_{1}(K_k)\stackrel{f
_{k*}}\longrightarrow H_{1}(K)\stackrel{\partial}\longrightarrow \\
\stackrel{\partial}\longrightarrow H_{0}(C_{\cdot}(K_k/\ker
(h_{k*}-id))\stackrel{t^k-1}\longrightarrow H_{0}(K_k)\stackrel{f
_{k*}}\longrightarrow H_{0}(K)\rightarrow  0,
\end{array}
\end{equation}
since
\[ \begin{array}{l}
im[H_{1}(C_{\cdot}(K_k/\ker
(h_{k*}-id))\stackrel{t^k-1}\longrightarrow H_{1}(K_k)]=A_k(G), \\
H_{1}(K)\simeq \mathbb Z, \\
H_{0}(C_{\cdot}(K_k/\ker (h_{k*}-id))\simeq \mathbb Z/k\mathbb Z, \\
H_{0}(K_k)\stackrel{f_{k*}}\simeq  H_{0}(K)\simeq \mathbb Z,
\end{array}
\]
are $\Lambda$-modules with trivial action of $t$ and exact
sequence (\ref{chain2}) is a sequence of $\Lambda$-homomorphisms
of $\Lambda$-modules.
\qed  \\

Now Theorem \ref{kmain1} follows from Lemmas \ref{splits} and
\ref{image}, since $\ker i_{k*}$ is generated by $[\bar \gamma]\in
H_1(X')_1$ and $f_{k*}([\bar\gamma])=k[\gamma]$.

Similarly, in Case II, we have $\ker i_{k*}=H_1(X'_k)_1$. Indeed,
$\ker i_{k*}$ is generated by $\bar\gamma $ and
$\bar\gamma_{\infty}\in H_1(X'_k)_1\simeq \mathbb Z$ and
$f_{k*}([\bar\gamma])=k[\gamma]$. Therefore $H_1(X'_k)_1$ is
generated by $[\bar\gamma]$.

As a consequence, we obtain that the restriction of $i_{k*}$ to
the submodule $A_k(G)$ of $H_1(X'_k,\mathbb Z)$ is an isomorphism
of $A_k(G)$ with $H_1(X_k,\mathbb Z)$. Therefore the following
lemma implies Theorem \ref{hmain2}.

\begin{lem} \label{isom} {\rm (}\cite{Gr-Ku}{\rm )}
The homomorphism $j_{k*} :H_1(X_k,\mathbb Q)\to H_1(\overline
X_k,\mathbb Q)$ is an isomorphism.
\end{lem}

\subsection{Corollaries of Theorems \ref{kmain1} and \ref{hmain2}}

\begin{cor} \label{ckmain1} Let $V$ be a knotted $n$-manifold,
$n\geq 1$, and ${f}_k:X_k\to S^{n+2}$ the cyclic covering branched
along $V$, $\deg f_k=k$. Then
\begin{itemize}
\item[$(i)$] the first Betti number $b_1(X_k)$ of $X_k$ is an even
number; \item[$(ii)$] if $k=p^r$, where $p$ is prime, then
$H_1(X_k,\mathbb Z)$ is finite; \item[$(iii)$] a finitely
generated abelian group $G$ can be realized as $H_1(X_k,\mathbb
Z)$ for some knotted $n$-manifold $V$, $n\geq 2$, if and only if
there is an automorphism $h\in \text{Aut}(G)$ such that
$h^k=\text{Id}$ and $h-\text{Id}$ is also an automorphism of $G$;
in particular, $H_1(X_2,\mathbb Z)$ is a finite abelian group of
odd order and any finite abelian group $G$ of odd order can be
realized as $H_1(X_2,\mathbb Z)$ for some knotted $n$-sphere,
$n\geq 2$.
\end{itemize}
\end{cor}

\proof It follows from Theorems \ref{main1}, \ref{kmain1},
\ref{even}, Propositions \ref{p^r-1}, Corollary \ref{sfera}, and
Examples \ref{ex1}, \ref{ex5}. \qed \\

Corollary \ref{finger} follows from Theorems \ref{main2} and
\ref{zfinite}.

Corollary \ref{ckmain2} is a simple consequence of Lemma
\ref{isom} and the following corollary, since the homomorphism
$j_{k*} :H_1(X_k,\mathbb Z)\to H_1(\overline X_k,\mathbb Z)$ is an
epimorphism and $H_1(\overline X_k,\mathbb Q)\simeq A_k(H)\otimes
\mathbb Q$.

\begin{cor} \label{ckmain3} Let $H$ be an algebraic
{\rm (}resp, Hurwitz or pseudo-holomorphic{\rm )} irreducible
curve in $\mathbb C\mathbb P^2$, $\deg H=m$, and ${\overline
f}_k:{\overline X}_k\to \mathbb C\mathbb P^2$ be a resolution of
singularities of the cyclic covering of degree $k$ branched along
$H$ and, maybe, alone a line "at infinity" $L$, and let
$X_k=\overline X_k\setminus E$. Then
\begin{itemize}
\item[$(i)$] the sequence of groups $$H_1(X_1,\mathbb Z), \dots,
H_1(X_k,\mathbb Z),\dots$$ has period $m$, that is,
$H_1(X_k,\mathbb Z) \simeq H_1(X_{k+m},\mathbb Z)$; \item[$(ii)$]
the first Betti number  $b_1(\overline X_k)=r_{k,\neq 1}$, where
$r_{k,\neq 1}$ is the number of roots of the Alexander polynomial
$\Delta (t)$ of the curve $H$ which are $k$-th roots of unity not
equal to $1$, in particular, $b_1(\overline X_k)$ is an even
number; \item[$(iii)$] if $k=p^r$, where $p$ is prime, then
$H_1(X_k,\mathbb Z)$ and $H_1(\overline X_k,\mathbb Z)$ are finite
groups; \item[$(iv)$] if $k$ and $m$ are coprime, then
$H_1(\overline X_k,\mathbb Z)=0$; \item[$(v)$] a finitely
generated abelian group $G$ can be realized as $H_1(X_k,\mathbb
Z)$ for some Hurwitz (resp., pseudo-holomorphic) curve $H$ if and
only if there is an automorphism $h\in \text{Aut}(G)$ such that
$h^d=\text{Id}$ and $h-\text{Id}$ is also an automorphism of $G$,
where $d$ is a divisor of $k$, and, moreover, if $G$ is realized
as $H_1(X_k,\mathbb Z)$ for a curve $H$, then $d$ is a divisor of
$\deg H$; in particular, $H_1(\overline X_2,\mathbb Z)$ is a
finite abelian group of odd order and any finite abelian group $G$
of odd order can be realized as $H_1(X_2,\mathbb Z)$ for some
Hurwitz (resp., pseudo-holomorphic) curve $H$ of even degree.
\end{itemize}
\end{cor}

\proof It follows from Theorems \ref{main2}, \ref{hmain2},
\ref{even}, \ref{period} and  Propositions \ref{p^r-1},
\ref{realization}. \qed \\

Note that there are plane algebraic curves $H$ for which the
homomorphisms $j_{k*} :H_1(X_k,\mathbb Z)\to H_1(\overline
X_k,\mathbb Z)$ are not isomorphisms.
\begin{ex} \label{exK3}
Let $H\subset \mathbb C\mathbb P^2$ be a curve of degree $6$ given
by equation $$Q^3(z_0,z_1,z_2)C^2(z_0,z_1,z_2)=0,$$ where $Q$ and
$C$ are homogeneous forms of $\deg Q=2$, $\deg C=3$ and the conic
and cubic, given by equations  $Q=0$ and $C=0$, meet transversally
at $6$ points. Then $A_2(H)\simeq \mathbb Z/3\mathbb Z$, but
$H_1(\overline X_2,\mathbb Z)=0$.
\end{ex}
\proof It is known (see \cite{Z}) that $\pi_1(\mathbb C\mathbb
P^2\setminus (H\cup L))\simeq \text{Br}_3$ as a $C$-group.
Therefore $A_2(H)\simeq \mathbb Z/3\mathbb Z$ (see Examples
\ref{ex2} and \ref{ex4}).

It is also well known that the minimal resolution of singularities
of two-sheeted covering of $\mathbb C\mathbb P^2$ branched along
$H$ is a $K3$-surface which is simply connected. \qed \\

Note also that in the case of knotted $n$-manifold $V\subset
S^{n+2}$ the sequence of homology groups $H_1(X_k,\mathbb Z)$,
$k\in \mathbb N$, is not necessary to be periodic. For example, if
$S^2\subset S^4$ is a knotted sphere for which $\pi_1(S^4\setminus
S^2)\simeq G_m$, where $G_m$ is a group considered in Example
\ref{ex5} (by Corollary \ref{sfera}, this group can be realized as
a group of knotted sphere), then $H_1(X_k,\mathbb Z)$ is the
cyclic group of order $(m+1)^k-m^k$ (see Example \ref{ex1}).

 \ifx\undefined\bysame
\newcommand{\bysame}{\leavevmode\hbox to3em{\hrulefill}\,}
\fi

\end{document}